\documentclass[12pt]{article}
\usepackage{amsmath,amsthm,amssymb}
\usepackage{graphicx}

\newcommand{\C}{\mbox{\rm \,l\kern-0.52em C}}
\newcommand{\Ce}{\rm \,l\kern-0.35em C}

\newcommand{\Inf}{{\rm Inf}}

\newtheorem{theorem}{Theorem}[section]
\newtheorem{deflem}[theorem]{Definition and Lemma}
\newtheorem{definition}[theorem]{Definition}
\newtheorem{prop}[theorem]{Proposition}

\newtheorem{cor}[theorem]{Corollary}

\newtheorem{lemma}[theorem]{Lemma}
\newtheorem{remark}[theorem]{Remark}

\renewenvironment{proof}{{\bf Proof:}}{\mbox{}\hfill $\Box$}

\theoremstyle{definition}
\date{}

\title{Local cyclic homology of group Banach algebras of ``non-positively curved'' discrete groups}

\author{Michael Puschnigg}
\begin{document}

\maketitle
\abstract{We calculate the local cyclic homology of group Banach-algebras of discrete groups acting properly, isometrically and cocompactly on a CAT(0)-space.}
\section{Introduction}

Cyclic homology of Banach algebras was introduced by Connes \cite{Co} as a homology theory which should allow, via the Chern-character, an approximate calculation of Banach $K$-theory by means of classical homological algebra.
While the theory and its formal properties were amply and quite successfully studied during the eighties and nineties of the last century, explicit calculations in cases not accessible by formal arguments
remained rare.\\
Let us give a heuristic explanation why the calculation of cyclic homology theories for Banach algebras turned out to be so difficult.
On the one hand Connes showed that cyclic homology of abstract algebras may be viewed as a classical derived functor \cite{Co}, and Nistor used this point of view to calculate the cyclic homology of group rings and algebraic crossed products \cite{Ni}. On the other hand Bott's Periodicity Theorem in its most general form, due to Cuntz and Higson \cite{CMR}, states that every additive, stable (under Morita equivalence), and split exact homotopy functor on the category of Banach-algebras is two-periodic under suspension. A stable and homotopy invariant homology theory for Banach algebras like local cyclic homology therefore cannot be a derived functor.
We view this as the reason why any attempt to extend Nistor's approach to analytic and local cyclic homology of group Banach algebras notoriously leads to serious convergence problems.\\
Therefore a new algebraic framework has to be set up in order to overcome these difficulties. It is provided by the approach to the periodic cyclic homology of group rings in \cite{Pu1}, which is based exclusively on ${\mathbb Z}/2{\mathbb Z}$-graded chain complexes. We extend it here to a topological setting. We finally introduce a notion of ``non-commutative'' derived functor and use it to obtain after a rather simple calculation our main result:
\begin{theorem} 
 Let $\Gamma$ be a discrete group which either acts properly, isometrically and cocompactly on a $CAT(0)$-space or which is word-hyperbolic.Then
 $$
 HC_*^{loc}(\ell^1(\Gamma))\,\simeq\,H_*(\Gamma,{\mathbb C}\Gamma_{tors})
 \eqno(1.1)
 $$
 where ${\mathbb C}\Gamma_{tors}$ denotes the submodule of ${\mathbb C}\Gamma$ spanned by the torsion elements, equipped with the adjoint action. 
\end{theorem}
In fact the analytic cyclic bicomplex of $\ell^1(\Gamma)$ becomes isomorphic to the finite dimensional vector space $H_*(\Gamma,{\mathbb C}\Gamma_{tors})$ in the derived ind-category \cite{Pu2}.\\
\\
Recall that the cyclic bicomplex of the group ring of a discrete group $\Gamma$ decomposes canonically into a direct sum of bicomplexes, 
labeled by the conjugacy classes of $\Gamma$. The isomorphism (1.1) is compatible with this decomposition.\\
\\
The right hand side of formula (1.1) is a classical derived functor. Its homogeneous part, i.e. the contribution of the conjugacy class of the unit, is obtained by taking the homology of the complex of $\Gamma$-coinvariants of 
a projective resolution $C_*(\Gamma,{\mathbb C})$(in our case the standard or Bar-resolution) of the constant $\Gamma$-module ${\mathbb C}$. 
We will identify the left hand side of formula (1.1) with a derived functor in a ``non-commutative'' sense. \\
\\
For us a ``non-commutative'' resolution of a $\Gamma$-algebra $A$ is a surjective homomorphism $R\to A$ of $\Gamma$-algebras with nilpotent kernel and such that $R$ is projective as $\Gamma$-module. By Goodwillie's Theorem or the homotopy invariance of periodic/analytic cyclic homology the induced morphism $CC_*(R)\,\to\,CC_*(A)$ of cyclic bicomplexes will be a chain-homotopy equivalence and $CC_*(R)$ will be a projective $\Gamma$-module. For $A={\mathbb C}$ we choose as resolution the algebra $R={\mathbb C}\langle\Gamma\rangle,$  given by the group ring ${\mathbb C}\Gamma$, equipped with the multiplication
$$
\begin{array}{lr}
a*b\,=\,\epsilon(a)\cdot b,  & \forall a,b \in {\mathbb C}\langle\Gamma\rangle,\\
\end{array}
\eqno(1.2)
$$
where $\epsilon:{\mathbb C}\Gamma\to{\mathbb C}$ is the trivial representation.
A suitable completion of the complex of $\Gamma$-coinvariants of $CC_*({\mathbb C}\langle\Gamma\rangle)$
equals then the homogeneous part of the analytic cyclic bicomplex of the group Banach algebra $\ell^1(\Gamma)$ and shows that the latter may be interpreted as a``derived functor''. In order to establish the main theorem it suffices therefore to compare the two resolutions. We show that the augmentation morphisms 
$$
CC_*({\mathbb C}\langle\Gamma\rangle)\,\longleftarrow\,CC_*({\mathbb C}\langle\Gamma\rangle)\otimes C_*(\Gamma,{\mathbb C})\,\longrightarrow\,C_*(\Gamma,{\mathbb C})
\eqno(1.3)
$$ 
extend to chain-homotopy equivalences of ind-Fr\'echet complexes of $\Gamma$-modules. Passing to $\Gamma$-coinvariants and taking homology yields our 
main result for the homogeneous parts.
A similar reasoning, following \cite{Pu1}, takes care of the contributions of the finitely many elliptic conjugacy classes of $\Gamma$. 
 The boundedness of the involved operators, which allows to pass to group Banach algebras, follows easily
  from the convexity of the distance function on a $CAT(0)$-space and the cocompactness of the $\Gamma$-action on it.\\
 \\
 To show the vanishing of the contributions of the hyperbolic conjugacy classes, we use contracting homotopy operators suggested by the work of Nistor \cite{Ni}. Their construction and boundedness relies on two facts.
 First, the set $Min(v)$ of points of minimal displacement under a hyperbolic isometry $v\in Isom(X)$ of a $CAT(0)$-space $X$ splits as a metric product of the real line (translated by $v$) and a space fixed pointwise by $v$.
Second, the action on $Min(v)$ of the centralizer of $v$ in $Isom(X)$ respects this splitting. As $\Gamma$ (acting cocompactly) does not contain parabolic isometries, the proof of the main result is complete.\\
\\
Our arguments readily apply to hyperbolic groups which we treat as well.\\
\\
We want to point out that our method does not provide any distinguished chain map of complexes inducing the isomorphism (1.1).  
The higher index theorem of Connes-Moscovici \cite{CM} implies however that (1.1) is compatible with the assembly map in Banach $K$-theory under the Chern-character. 
Therefore one may deduce from Lafforgue's work \cite{La} and our result that the Chern-character gives rise to an isomorphism 
$$
ch:\,K_*(\ell^1(\Gamma))\otimes_{\mathbb Z}{\mathbb C}\,\overset{\simeq}{\longrightarrow}\,HC_*^{loc}(\ell^1(\Gamma))
\eqno(1.4)
$$
for the class of discrete groups considered here.\\
\\
The content of the paper is as follows. Section 2 recalls the Bar- and Rips-complexes attached to a discrete group and discusses some contracting homotopies of these resolutions.
Section 3 introduces ``standard'' algebras and uses them to construct the cyclic or ``non-commutative'' resolutions of the constant $\Gamma$-algebra $\mathbb C$. In section 4 various locally convex topologies 
on these resolutions are introduced and it is shown that the homotopy equivalences in (1.3) are continuous with respect to these topologies. In section 5 we deduce our main theorems from the results of the previous sections.
 In a final remark we indicate a correction of an error committed in \cite{Pu4}.

\section{The standard resolution}

\begin{definition}
The {\bf Bar-complex} $\Delta_\bullet(X)$ of a set $X$ is the simplicial set 
with $n$-simplices $\Delta_n(X)\,=\,X^{n+1}$, face maps
$$
\partial_i([x_0,\ldots,x_n])\,=\,[x_0,\ldots,x_{i-1},x_{i+1},\ldots,x_n],
\eqno(2.1)
$$
and degeneracy maps
$$
s_j([x_0,\ldots,x_n])=[x_0,\ldots,x_j,x_j,\ldots,x_n]\,\,\text{for}\,\,0\leq i,j\leq n\in{\mathbb N}.
\eqno(2.2)
$$
The {\bf support} of a Bar-simplex is
$
Supp([x_0,\ldots,x_n])\,=\,\{x_0,\ldots,x_n\}\subset X.
$
\end{definition}

\begin{definition}
Let $(X,d)$ be a metric space and let $R\geq 0$. The {\bf Rips-complex} $\Delta_\bullet^R(X)$
of $(X,d)$ is the simplicial subset of the Bar-complex $\Delta_\bullet(X)$ given by the Bar-simplices of diameter at most $R:$
$$
\Delta_n^R(X)\,=\,\{[x_0,\ldots,x_n]\in \Delta_n(X),\,d(x_i,x_j)\leq R,\,0\leq i,j\leq n\}.
\eqno(2.3)
$$
\end{definition}

Every map of sets $f:X\to Y$ gives rise to a simplicial map 
$$
f_\bullet:\,\Delta_\bullet(X)\to\Delta_\bullet(Y),\,[x_0,\ldots,x_n]\mapsto [f(x_0),\ldots,f(x_n)].
\eqno(2.4)
$$
In particular, every group action on the set $X$ gives rise to a simplicial action on the Bar-complex $\Delta_\bullet(X)$ and every isometric group action on a metric space $(X,d)$ gives rise 
to a simplicial action on the Rips-complexes $\Delta_\bullet^R(X),\,R>0$.\\

\begin{definition}
The {\bf Bar chain complex} $C_*(X,{\mathbb C})$ of a set $X$ is given by the complex vector space with basis $\Delta_*(X)$. Its differentials are given by the alternating sum of the linear operators induced by the face maps. It is augmented by 
$$
\epsilon:\,C_*(X,{\mathbb C})\,\to\,{\mathbb C},\,[x]\mapsto 1,\,x\in X.
\eqno(2.5)
$$
The {\bf Rips chain complexes} $C_*^R(X,{\mathbb C})$ of a metric space are defined similarly. They are subcomplexes of the Bar chain complex.
\end{definition}

\begin{lemma}
Let $\varphi_*,\psi_*:C_*(X,{\mathbb C})\to C_*(Y,{\mathbb C})$ be chain maps of Bar complexes such that  $\epsilon_Y\circ\varphi=\epsilon_Y\circ\psi$, where $\epsilon_Y:C_*(Y,{\mathbb C})\to{\mathbb C}$ denotes the augmentation. Then the linear operator
$$
\begin{array}{cccc}
h_{Bar}(\varphi,\psi): & C_*(X,{\mathbb C}) & \to & C_{*+1}(Y,{\mathbb C}) \\
& [x_0,\ldots, x_n] & \mapsto & \underset{i=0}{\overset{n}{\sum}}\,(-1)^i\,[\varphi_i(x_0,\ldots,x_i),
\psi_{n-i}(x_i,\ldots,x_n)] \\
\end{array}
\eqno(2.6)
$$
defines a natural chain homotopy between $\varphi$ and $\psi$:
$$
\psi_*-\varphi_*\,=\,\partial\circ h_{Bar}(\varphi,\psi)\,+\,h_{Bar}(\varphi,\psi)\circ\partial.
\eqno(2.7)
$$
In particular, if $G$ is a group acting on $X$ and $Y$, and if $\varphi_*$ and $\psi_*$ are $G$-equivariant, then $h(\varphi,\psi)$ is $G$-equivariant as well.
\end{lemma}

The lemma implies that the augmentation map $\epsilon_X: C_*(X,{\mathbb C}) \to {\mathbb C}$ is a chain homotopy equivalence.

\begin{lemma}
The antisymmetrization operator
$$
\begin{array}{cccc}
\pi_{alt}: & C_*(X,{\mathbb C}) & \to & C_*(X,{\mathbb C}) \\
& & & \\
& [x_0,\ldots,x_n] & \mapsto & \frac{1}{(n+1)!}\underset{\sigma\in\Sigma_{n+1}}{\sum}\,(-1)^{\epsilon(\sigma)}\,[x_{\sigma(0)},\ldots,x_{\sigma(n)}] \\
\end{array}
\eqno(2.8)
$$
is a chain map which preserves the Rips subcomplexes and equals the identity in degree zero. In particular it is canonically chain homotopic to the identity by the previous lemma.
\end{lemma}

\begin{definition}
\begin{itemize}
\item[a)]
The {\bf alternating Bar chain complex} $C_*(X,{\mathbb C})^{alt}$ of $X$ is the image 
of the full Bar chain complex under the antisymmetrization operator. 
\item[b)]
Similarly the {\bf alternating Rips chain complexes} $C_*^{R}(X,{\mathbb C})^{alt},\,R>0,$ are the images of the corresponding full Rips chain complexes under the antisymmetrization operator.
\end{itemize}
\end{definition}

The alternating Bar- and Rips-chain complexes are at the same time quotients and subcomplexes (deformation retracts) of the full Bar- and Rips chain complexes.\\
\\
Let $\varphi:\,C_0(X,{\mathbb C})\to C_0(Y,{\mathbb C})$ be a linear map which commutes with augmentations, i.e. $\epsilon_Y\circ\varphi_0=\epsilon_X$. Then $\varphi_0$ extends to a chain map 
of Bar chain-complexes given by 
$$
\varphi_n\,=\,\varphi_0^{\otimes^{n+1}}:\,C_n(X,{\mathbb C})\simeq C_0(X,{\mathbb C})^{\otimes^{n+1}} \to C_0(Y,{\mathbb C})^{\otimes^{n+1}}\simeq C_n(Y,{\mathbb C}),\,n\geq 0.
\eqno(2.9)
$$ 
\\
In the sequel we will deal with two classes of discrete groups \cite{BH}:
\begin{itemize}
\item[(A)] Groups acting properly, isometrically and cocompactly on a\\ $CAT(0)$-space.
\item[(B)] Word-hyperbolic groups.
\end{itemize}
Case A:\\Let $\Gamma$ be a discrete group acting properly, isometrically and cocompactly on the $CAT(0)$-space $X$. This implies that $X$ is proper, i.e. every closed bounded subset of $X$ is compact. Fix once and for all a relatively compact fundamental domain $D$ of the $\Gamma$-action on $X$, i.e. a relatively compact Borel subset $D\subset X$ such that every $\Gamma$-orbit of $X$ intersects $D$ in a single point. We choose a base point $x\in D$ of $X$ and equip $\Gamma$ with the proper, left-$\Gamma$-invariant pseudo-metric 
$$
d_\Gamma:\Gamma\times\Gamma\to{\mathbb R}_+,\,d_\Gamma(g,h)=d_X(gx,hx).
\eqno(2.10)
$$
It is a metric if the (finite) stabilizer $\Gamma_x$ of $x$ is trivial.\\
\\
Let $(X,d)$ be a complete $CAT(0)$-space. Let $Y\subset X$ be a closed, convex subset. Let $$\pi'_{Y}:X\to Y$$  be the canonical projection which sends a point $x'\in X$ to the unique point of $Y$ at minimal distance from $x$.
The projection does not increase distances, i.e.
$$
d(\pi'_{Y}(x'),\pi'_{Y}(x''))\,\leq\,d(x',x''),\,\,\forall\,x',x''\in X.
\eqno(2.11)
$$
For $v\in\Gamma_x$ let $X^v\subset X$ be the set of $v$-fixed points in $X$. It is closed, convex and not empty \cite{BH} II.2.8. 
Denote by $C_*(\Gamma,v)$ the subcomplex of $C_*(\Gamma,{\mathbb C})$ spanned by 
$$
\Delta_\bullet(\Gamma,v)\,=\,\{\alpha\in\Delta_{\bullet}(\Gamma),\,gD\cap X^v\neq\emptyset,\,\forall g\in Supp(\alpha)\}
\eqno(2.12)
$$
and put 
$C^R_*(\Gamma,v)=C_*(\Gamma,v)\cap C_*^R(\Gamma,{\mathbb C})$.
For $y\in X^v$ and $n\in{\mathbb N}$ let 
$$
\pi_{(y,v,n)}: C_0(\Gamma,{\mathbb C})\to C_0(\Gamma,{\mathbb C})
\eqno(2.13)
$$
be the linear map characterized by 
$$
\pi_{(y,v,n)}([g])\,=\,\frac{1}{\vert S_{(y,v,n)}(g)\vert}\cdot\underset{h\in S_{(y,v,n)}(g)}{\sum}\,[h]
\eqno(2.14)
$$
for all $g\in\Gamma$, where 
$$
S_{(y,v,n)}(g)\,=\,\{h\in\Gamma,\,\pi'_{X^v}\circ\pi'_{B(y,n)}(gx)\in hD\}.
\eqno(2.15)
$$
It is equivariant with respect to the action of the finite cyclic subgroup $U\subset\Gamma$ spanned by $v\in\Gamma$ because every ball $B(y,r)$centered at $y$ and the fixed point set
$X^v$ are invariant under $U$.
As it commutes with augmentations it gives rise to a $U$-equivariant chain map
$$
(\pi_{(y,v,n)})_*:\,C_*(\Gamma,v)\,\longrightarrow C_*(\Gamma,v).
\eqno(2.16).
$$
Because ``orthogonal'' projections do not increase distances \cite{BH} II.2.5 one has
$$
d_\Gamma(h,h')\,\leq\,d_\Gamma(g,g')+\vert n-n'\vert+2\,diam(D)
\eqno(2.17)
$$
for all $h\in S_{(y,v,n)}(g),\,h'\in S_{(y,v,n')}(g'),\,g,g'\in\Gamma$. Therefore (2.16) restricts to 
$$
(\pi_{(y,v,n)})_*:\,C^R_*(\Gamma,v)\,\longrightarrow C^{R+2\,diam(D)}_*(\Gamma,v),\,\,\,\forall R>0.
\eqno(2.18)
$$

\begin{lemma}
The chain map
$$
(\pi_{(y,v)})_*:\,C_*(\Gamma,v)\,\longrightarrow C_*(\Gamma,v)
\eqno(2.19)
$$ 
induced by the linear projection 
$$
\begin{array}{cccc}
\pi_{(y,v)}: & C_0(\Gamma,{\mathbb C}) & \to & C_0(\Gamma,v) \\
 & & & \\
 & [g] & \mapsto & \frac{1}{\vert\Gamma_y\vert}\underset{h\in\Gamma_y}{\sum}\,[h\widetilde{h}], \\
\end{array}
\eqno(2.20)
$$
where $\widetilde{h}\in\Gamma$ is any element satisfying $y\in\widetilde{h}D$, 
is $U$-equivariantly chain-homotopic to the identity.
The same holds for alternating Bar chain complexes.
\end{lemma}

\begin{proof}
In the notations of (2.6) and (2.16) put
$$
h_{(y,v)}\,=\,\underset{n=0}{\overset{\infty}{\sum}}\,h_{Bar}(\pi_{(y,v,n)},\pi_{(y,v,n+1)}):\,C_*(\Gamma,v)\,\longrightarrow\,C_{*+1}(\Gamma,v).
\eqno(2.21)
$$
Note that $h_{Bar}(\pi_{(y,v,r)},\pi_{(y,v,r+1)})(\alpha)=0$ for Bar-simplices $\alpha\in\Delta_\bullet(\Gamma,v)$ supported in $\overline{B(y,r)}$, so that this operator is well defined. It 
 defines a $U$-equivariant chain homotopy between $(\pi_{(y,v)})_*=(\pi_{y,v,0})_*$ and the identity. Using the homotopy operator
 $$
 h^{alt}_{(y,v)}\,=\,\pi_{alt}\circ h_{(y,v)}:\,C_*(\Gamma,v)^{alt}\,\longrightarrow\,C_{*+1}(\Gamma,v)^{alt},
 \eqno(2.22)
 $$
 one gets the corresponding assertion for alternating Bar-complexes.
\end{proof}

The homotopy operator (2.21) restricts for all $R>0$ to a linear map
$$
h_{(y,v)}:\,C^R_*(\Gamma,v)\,\longrightarrow\,C^{R+2\,diam(D)+1}_{*+1}(\Gamma,v).
\eqno(2.23)
$$
Case B:\\
Let $(\Gamma,S)$ now be a $\delta$-hyperbolic group, i.e. $\Gamma$ is a group with finite symmetric set $S$ of generators such that the Cayley-graph ${\mathcal G}(\Gamma,S)$ is a $\delta$-hyperbolic geodesic metric space \cite{BH}. Hyperbolicity does not depend on the choice of $S$ (but the hyperbolicity constant $\delta\geq 0$ does). For a subset $Y\subset X$ let
$$
\begin{array}{cccc}
\pi_{Y}: & C_0(\Gamma,{\mathbb C}) & \to & C_0(\Gamma,{\mathbb C}) \\
 & & & \\
 & [g] & \mapsto & \frac{1}{\vert S'_{Y}(g)\vert}\underset{h\in S'_{Y}(g)}{\sum}\,[h] \\
\end{array}
\eqno(2.24)
$$
where $S'_{Y}(g)=\{h\in Y, d(g,h)=d(g,Y)\}$. One has
$$
d_S(h,h')\,\leq\,d_S(g,g')+\vert n-n'\vert+2\delta
\eqno(2.25)
$$
for all $h\in S'_{B(y,n)}(g),\,h'\in S'_{B(y,n')}(g'),\,g,g'\in\Gamma$, where $B(y,r)$ denotes the closed $r$-ball around $y\in Z_v$. Therefore 
$$
(\pi_{B(y,n)})_*:\,C^R_*(\Gamma,{\mathbb C})\,\longrightarrow C^{R+2\delta}_*(\Gamma,{\mathbb C}),\,\,\,\forall R\geq 0.
\eqno(2.26)
$$
The linear operator
$$
h_{(y,v)}\,=\,\underset{n=0}{\overset{\infty}{\sum}}\,h_{Bar}(\pi_{B(y,n)},\pi_{B(y,n+1)}):\,C_*(\Gamma,{\mathbb C})\,\longrightarrow\,C_{*+1}(\Gamma,{\mathbb C}).
\eqno(2.27)
$$
is a chain homotopy between the identity and the constant map onto $\{y\}$.
It restricts for every $R>0$ to a linear map
$$
h_{(y,v)}:C^R_*(\Gamma,{\mathbb C})\,\longrightarrow\,C^{R+2\delta+1}_{*+1}(\Gamma,{\mathbb C}).
\eqno(2.28)
$$
In all cases we have for $g,v\in\Gamma,\,\vert v\vert<\infty, y\in X^v$ or $y\in Z_v$  the identities
$$
\begin{array}{ccc}
g\circ \pi_{(y,v)}\circ g^{-1}\,=\,\pi_{(gy,gvg^{-1})}, &  g\circ h_{(y,v)}\circ g^{-1}\,=\,h_{(gy,gvg^{-1})}, & g\circ h^{alt}_{(y,v)}\circ g^{-1}\,=\,h^{alt}_{(gy,gvg^{-1})}.\\
\end{array}
\eqno(2.29)
$$

\section{The cyclic resolution}

\subsection{The cyclic bicomplex}

Let $R$ be a unital and associative complex algebra. An $R$-algebra is a unitary\\ $R$-bimodule $A$, equipped with an $R$-bimodule homomorphism $m:\,A\otimes_R A\,\to\,A,$
 turning $A$ into a (not necessarily unital) complex algebra \cite{Pu1}. A homomorphism of $R$-algebras is a homomorphism of the underlying complex algebras, which is at the same time a homomorphism of $R$-bimodules.\\
 \\
Recall that the graded $R$-bimodule of  {\bf algebraic differential forms} of $A$ over $R$ is defined as
$$
\begin{array}{ccc}
\Omega^n(A:R) & \simeq & A^{\otimes_R^{n+1}} \oplus A^{\otimes_R^{n}} \\
& &  \\
a_0da_1\ldots da_n & \leftrightarrow & a_0\otimes a_1\otimes\ldots\otimes a_n  \\
& &  \\
da_1\ldots da_n & \leftrightarrow &  a_1\otimes\ldots\otimes a_n  \\
\end{array}
\eqno(3.1)
$$
with $R$-bimodule-structure induced by the natural one on the tensor powers of $A$. 
The commutator quotient with respect to this bimodule structure is the graded vector space 
$$
\Omega^*(A:R)_\natural\,=\,\Omega^*(A:R)/[\Omega^*(A:R),R].
\eqno(3.2)
$$
Under the canonical projection 
$$
\Omega^*A\,\to\,\Omega^*(A:R)_\natural
\eqno(3.3)
$$
the well known Hochschild- and Connes-operators 
$$
\begin{array}{cc}
b:\,\Omega^*A\,\to\,\Omega^{*-1}A, & B:\,\Omega^*A\,\to\,\Omega^{*+1}A
\end{array}
\eqno(3.4)
$$
descend to linear operators 
$$
\begin{array}{cc}
b:\,\Omega^*(A:R)_\natural\,\to\,\Omega^{*-1}(A:R)_\natural, & B:\,\Omega^*(A:R)_\natural\,\to\,\Omega^{*+1}(A:R)_\natural,
\end{array}
\eqno(3.5)
$$
satisfying the usual identities
$$
b^2\,=\,B^2\,=\,bB+Bb\,=\,0.
\eqno(3.6)
$$
The ${\mathbb Z}/2{\mathbb Z}$-graded {\bf cyclic bicomplex } of $A$ over $R$ is given by 
$$
CC_*(A:R)\,=\,\left(\Omega^*(A:R)_\natural,\,\,b+B\right).
\eqno(3.7)
$$
It is graded by the parity of differential forms. The canonical projection (3.3) gives rise to an epimorphism
$$
\pi:\,CC_*(A)\,\longrightarrow\,CC_*(A:R)
\eqno(3.8)
$$
of cyclic bicomplexes.\\

The chain maps of cyclic bicomplexes induced by two smoothly homotopic $R$-algebra
homomorphisms are naturally chain homotopic by a well known homotopy formula which we recall now \cite{Lo}, 4.1.8. 

\begin{prop}
Let $F_t:A\to B,\,0\leq t\leq 1,$ be a family of $R$-algebra homomorphisms depending smoothly on a parameter
$t$. Then the Cartan homotopy formula
$$
\frac{\partial}{\partial t}\,CC_*(F_t)\,=\,(b+B)(h_t+H_t)+(h_t+H_t)(b+B)
\eqno(3.9)
$$
holds for the linear operators 
$$
\begin{array}{ccc}
h_t:\Omega^*(A:R)_\natural\to\Omega^{*-1}(A:R)_\natural & \text{and} & H_t:\Omega^*(A:R)_\natural\to\Omega^{*+1}(A:R)_\natural
\end{array}
$$
given by the formulas
$$
h_t(a^0da^1\ldots da^n)\,=\,(-1)^{n-1}\overset{\bullet}{F_t}(a^n)F_t(a^0)dF_t(a^1)\ldots dF_t(a^{n-1})
\eqno(3.10)
$$
and
$$
H_t(a^0da^1\ldots da^n)\,=\,
$$
$$
\underset{1\leq i\leq j\leq n}{\sum}\,(-1)^{in+1}dF_t(a^i)dF_t(a^{i+1})\ldots dF_t(a^{j-1})
d\overset{\bullet}{F_t}(a^j)dF_t(a^{j+1})\ldots dF_t(a^n)dF_t(a^0)\ldots dF_t(a^{i-1}).
\eqno(3.11)
$$
\end{prop}

\subsection{Standard algebras}

Let $X$ be a set and let ${\mathbb C}\langle X\rangle$ be the complex vector space over $X$. The projection
$$
p:\,X\times X\to X,\,\,(x,y)\mapsto y
\eqno(3.12)
$$
turns ${\mathbb C}\langle X\rangle$ into an associative complex algebra with multiplication
$$
m:\,{\mathbb C}\langle X\rangle\otimes_{\mathbb C}{\mathbb C}\langle X\rangle\,=\,{\mathbb C}
\langle X\times X\rangle\,\overset{{\mathbb C}p}{\longrightarrow}\,{\mathbb C}\langle X\rangle.
\eqno(3.13)
$$
If $X$ consists of a single element the algebra thus obtained is canonically isomorphic to the complex field. Every map of sets $f:X\to Y$ gives rise to an algebra homomorphism 
$$
{\mathbb C}\langle f\rangle:{\mathbb C}\langle X\rangle\to {\mathbb C}\langle Y\rangle.
\eqno(3.14)
$$
In particular the constant map to a point induces an augmentation homomorphism
$$
\epsilon_X:{\mathbb C}\langle X\rangle\to{\mathbb C}
\eqno(3.15)
$$
The multiplication in ${\mathbb C}\langle X\rangle$ is characterized by the identity

$$
a\cdot b\,=\,\epsilon_X(a)b,\,\,\,\forall a,b\in{\mathbb C}\langle X\rangle.
\eqno(3.16)
$$

\begin{lemma}
\begin{itemize}
\item[a)]
Let $\varphi:{\mathbb C}\langle X\rangle
\to{\mathbb C}\langle Y\rangle$ be a linear map. Then $\varphi$ is a non-vanishing homomorphism of algebras iff it is compatible with augmentations, i.e. if
$$
\epsilon_Y\circ\varphi=\epsilon_X.
\eqno(3.17)
$$
\item[b)]
Any two non-vanishing algebra homomorphisms $\varphi,\psi:{\mathbb C}\langle X\rangle
\to{\mathbb C}\langle Y\rangle$ are smoothly homotopic via the affine homotopy 
$$
F_t=(1-t)\varphi+t\psi,\,0\leq t\leq 1
\eqno(3.18)
$$
of algebra homomorphisms.
\end{itemize}
\end{lemma}

\begin{proof}
Let $\varphi:{\mathbb C}\langle X\rangle
\to{\mathbb C}\langle Y\rangle$ be a linear map and let $x,x'\in X$. Then
$$
\begin{array}{ccc}
\varphi(x\cdot x')=\varphi(\epsilon_X(x)x')=\epsilon_X(x)\varphi(x') & \text{and} &
\varphi(x)\cdot\varphi(x')=\epsilon_Y(\varphi(x))\varphi(x').
\end{array}
$$

So on the one hand the condition $\epsilon_Y\circ\varphi=\epsilon_X$ implies that $\varphi$ is a homomorphism of algebas. Moreover it does not vanish because $\epsilon_X\neq 0$. If on the other hand $\varphi$ is a homomorphism of algebras and $x'\in X$ is chosen such that $\varphi(x')\neq 0$ 
the two equations above imply that $\epsilon_Y\circ\varphi=\epsilon_X$. The second assertion follows immediately from the first one by verifying the compatibility of the linear maps $F_t,\,0\leq t\leq 1,$ with augmentations.
\end{proof}

It follows that the augmentation map from ${\mathbb C}\langle X\rangle$ to the complex field is a smooth homotopy equivalence.  

\subsection{The cyclic resolution}

The integrated version of the Cartan homotopy formula 3.1 for the affine homotopy (3.18) between non-vanishing homomorphisms of standard algebras looks as follows.

\begin{lemma}
Let $f,g:{\mathbb C}\langle X\rangle\to {\mathbb C}\langle Y\rangle$ be non-vanishing algebra homomorphisms. Then 
$$
CC({\mathbb C}\langle g\rangle)-CC({\mathbb C}\langle f\rangle)\,=\,(b+B)\circ (h+H)(f,g)+
(h+H)(f,g)\circ(b+B)
\eqno(3.19)
$$
with $h(f,g):\Omega^*({\mathbb C}\langle X\rangle)\to\Omega^{*-1}({\mathbb C}\langle Y\rangle)$ and 
$H(f,g):\Omega^*({\mathbb C}\langle X\rangle)\to\Omega^{*+1}({\mathbb C}\langle Y\rangle)$ given for all $x_0,x_1,\ldots,x_n\in X$ and all $n\in{\mathbb N}$ by the formulas
$$
h(f,g)(dx_1\ldots dx_n)\,=\,\frac{(-1)^{n-1}}{n}\underset{k=0}{\overset{n-1}{\sum}}\,
\left(\begin{matrix}
n-1 \\
k\\
\end{matrix}\right)^{-1}\underset{\underset{\vert S\vert=k}{S\subset\{1,\ldots,n-1\}}}{\sum}\,
(g(x_n)-f(x_n))du_S(x_1)\ldots du_S(x_{n-1}),
\eqno(3.20)
$$
and
$$
H(f,g)(x_0dx_1\ldots dx_n)\,=\,\underset{1\leq i\leq j\leq n}{\sum}\frac{1}{n+1}\underset{k=0}{\overset{n}{\sum}}\,
\left(\begin{matrix}
n \\
k\\
\end{matrix}\right)^{-1}\underset{\underset{\vert S\vert=k}{S\subset\{0,\ldots,n\}\backslash\{j\}}}{\sum}\omega_{i,j,S}
\eqno(3.21)
$$
with
$$
\omega_{i,j,S}=(-1)^{ni+1}du_S(x_i)\ldots du_S(x_{j-1})d(g(x_j)-f(x_j))du_S(x_{j+1})\ldots
du_S(x_n)du_S(x_0)\ldots du_S(x_{i-1})
\eqno(3.22)
$$
where 
$$
u_S(x_j)\,=\,\begin{cases}
f(x_j) & j\in S, \\
g(x_j) & j\notin S. \\
\end{cases}
\eqno(3.23)
$$
Moreover
$$
h(f,g)(x_0dx_1\ldots dx_n)=0
\eqno(3.24)
$$
and
$$
H(f,g)(dx_1\ldots dx_n)=0.
\eqno(3.25)
$$
Here we identify elements of $X$ with the corresponding basis vectors in ${\mathbb C}\langle X\rangle$. 
\end{lemma}

This follows from a straightforward calculation by noting that for integers $0\leq k\leq n$
$$
(n+1)\left(\begin{matrix}
n \\
k\\
\end{matrix}\right)\underset{0}{\overset{1}{\int}}\,(1-t)^kt^{n-k}dt=1.
$$

In the sequel the notations of Section 2 will be understood. 

\begin{lemma}
Let $\Gamma$ be an isometry group of a $CAT(0)$-space as in Section 2 or let $\Gamma$ be word-hyperbolic. The endomorphism
of standard algebras induced by the constant map of $\Gamma$ onto $\{g'\}\subset\Gamma$ gives rise to a chain map 
$$
\pi_{g'}^{cyc}:\,CC_*({\mathbb C}\langle\Gamma\rangle)\,\longrightarrow\,CC_*({\mathbb C}\langle\Gamma\rangle)
$$
of cyclic bicomplexes which is chain homotopic to the identity.
\end{lemma}

\begin{proof}
For an isometry group of a $CAT(0)$-space put
$$
h^{cyc}_{g'}\,=\,(h+H)(\pi_{g'},\pi_{(g'x,e,0)})\,+\,\underset{n=0}{\overset{\infty}{\sum}}\,(h+H)(\pi_{(g'x,e,n)},\pi_{(g'x,e,n+1)})
\eqno(3.26)
$$
(see (2.14.)-(2.16)). This defines an operator
$h^{cyc}_{g'}:\,CC_*({\mathbb C}\langle\Gamma\rangle)\,\longrightarrow\,CC_{*+1}({\mathbb C}\langle\Gamma\rangle)$
because $h(\pi_{(g'x,e,r)},\pi_{(g'x,e,r+1)})(\beta)=H(\pi_{(g'x,e,r)},\pi_{(g'x,e,r+1)})(\beta)=0$ for $\beta\in\Omega^*({\mathbb C}\langle\Gamma\rangle)$ supported in $B(g'x,r)$, so that the formal sum (3.26) becomes finite when evaluated on a given cyclic simplex. It defines a chain homotopy between $\pi_{g'}^{cyc}$ and the identity. In the case of word-hyperbolic groups we use the homotopy operator 
$$
h^{cyc}_{g'}\,=\,\underset{n=0}{\overset{\infty}{\sum}}\,(h+H)(\pi_{(g',n)},\pi_{(g',n+1)})
\eqno(3.27)
$$
instead (the notations being those of (2.24)-(2.26)).
\end{proof}

 In all cases these operators are compatible with the $\Gamma$-action, i.e. 
$$
\begin{array}{cccc}
g\circ \pi_{g'}^{cyc}\circ g^{-1}\,=\,\pi_{gg'}^{cyc} & \text{and} & g\circ h^{cyc}_{g'}\circ g^{-1}\,=\,h_{gg'}^{cyc} & \text{for all}\,\,g,g'\in\Gamma.\\
\end{array}
\eqno(3.28)
$$

\section{The tilting complexes}

\subsection{Ind-complexes}

Recall that a category $I$ is filtrant if 

\begin{itemize}
\item For any two objects $i,j\in ob(I)$ there exists a diagram $i\to k\leftarrow j$ in $I$. \\
\item For any two morphisms $\alpha,\beta\in mor_I(i,j)$ there exists a morphism\\ $\gamma\in mor_I(j,k')$ such that $\gamma\circ\alpha=\gamma\circ\beta$.
\end{itemize}

Recall the ind-category $ind({\mathcal C})$  associated to a given category $\mathcal C$. It has the family of covariant functors from small filtrant categories to $\mathcal C$ as objects, and 
the set of morphisms between two ind-objects is given by
$$
Mor_{ind({\mathcal C})}("\underset{i\in I}{\lim}\,"X_i,"\underset{j\in J}{\lim}\,"Y_j)\,=\,\underset{\underset{i\in I}{\longleftarrow}}{\lim}\,\underset{\underset{j\in J}{\longrightarrow}}{\lim}\,Mor_{{\mathcal C}}(X_i,Y_j),
\eqno(4.1)
$$
where the limits on the right hand side are taken in the category of sets. Note that the ind-category admits small filtrant colimits.\\
\\
An ind-Fr\'echet-complex is an object in the ind-category of ${\mathbb Z}/2{\mathbb Z}$-graded chain complexes of complex Fr\'echet spaces.
A chain homotopy between morphisms $\varphi,\psi:"\underset{i\in I}{\lim}"C_*^{(i)}\,\to\,"\underset{j\in J}{\lim}" \widetilde{C}_*^{(j)}$ is an odd morphism $h:"\underset{i\in I}{\lim}"C_*^{(i)}\,\to
\,"\underset{j\in J}{\lim}"\widetilde{C}_{*+1}^{(j)}$ of ${\mathbb Z}/2{\mathbb Z}$-graded ind-Fr\'echet spaces satisfying $\psi-\varphi=\partial\circ h+h\circ\partial$. 

\subsection{Some locally convex topologies}

Let $\Gamma$ be a discrete group acting properly, isometrically and cocompactly on the $CAT(0)$-space $X$, or suppose that $\Gamma$ is word-hyperbolic. In the first case let $x\in X$ be a base point, and let $D\subset X$ be a relatively compact fundamental domain of the $\Gamma$-action containing $x$. In the second case let $d_\Gamma$ be the word-metric associated to a finite symmetric set of generators of $\Gamma$. Let $v\in\Gamma_x$ (resp. $v\in\Gamma$) be an element of finite order, let $Z_v\subset\Gamma$ be its centralizer, and let $U=v^{{\mathbb Z}}\subset\Gamma$ be the finite cyclic subgroup generated by $v$. The notations of sections 2 and 3 are understood.

\begin{deflem} 
\begin{itemize}
\item[a)]  Let $\Gamma$ be a discrete isometry group of a $CAT(0)$-space as in Section 2.
 For $R> 0$ let 
$$
{\mathcal C}^{R}_*(\Gamma,v)=\overline{C^{R}_*(\Gamma,v)}\subset\left(\ell^1(\Delta_\bullet(\Gamma)),\partial_{bar}\right)
\eqno(4.2)
$$ 
and denote by ${\mathcal C}^{R}_*(\Gamma,v)^{alt}$ its image under antisymmetrization.\\
\item[b)] For a word-hyperbolic group $\Gamma$ and $R>0$ put
$$
{\mathcal C}^{R}_*(\Gamma,v)=\overline{C^{R}_*(Z_v,d_\Gamma)}\subset\left(\ell^1(\Delta_\bullet(\Gamma)),\partial_{bar}\right)
\eqno(4.3)
$$ 
and denote by ${\mathcal C}^{R}_*(\Gamma,v)^{alt}$ its image under antisymmetrization.\\
\item [c)]Put
$$
{\mathcal C}_*(\Gamma,v)\,=\,"\underset{R\to \infty}{\lim}"\,{\mathcal C}^{R}_*(\Gamma,v),
\eqno(4.4)
$$
and 
$$
{\mathcal C}^{alt}_*(\Gamma,v)\,=\,"\underset{R\to \infty}{\lim}"\,{\mathcal C}^{R}_*(\Gamma,v)^{alt}.
\eqno(4.5)
$$
These are ind-Banach chain complexes of $Z_v$-modules.
\end{itemize}
\end{deflem}

\begin{deflem}
Let $\Gamma$ be a discrete isometry group of a $CAT(0)$-space as in Section 1 or let $\Gamma$ be word-hyperbolic. 
Let $v\in\Gamma$ be an element of finite order.
\begin{itemize}
\item[a)]
Denote by
$$
\Delta^{cyc}_n(\Gamma)\,=\,\Delta^{cyc}_{n}(\Gamma)'\cup\Delta^{cyc}_{n}(\Gamma)''
\eqno(4.6)
$$
$$
=\,\{\langle g_0\rangle d\langle g_1\rangle\ldots d\langle g_n\rangle,\,g_0,\ldots,g_n\in\Gamma\}\,
\cup\,\{ d\langle g_1\rangle\ldots d\langle g_n\rangle,g_1,\ldots,g_n\in\Gamma\}
$$
the canonical  basis of $\Omega^n({\mathbb C}\langle\Gamma\rangle),\,n\geq 0$. Its elements are called cyclic\\ $n$-simplices.
The support of cyclic simplices is the set of its vertices.
\item[b)]
The $v$-weight of cyclic simplices is given by 
$$
\vert\langle g_0\rangle d\langle g_1\rangle\ldots d\langle g_n\rangle\vert_v\,=\,d_\Gamma(g_0,g_1)+\ldots+d_\Gamma(g_{n-1},g_n)+d_\Gamma(g_n,vg_0)
\eqno(4.7)
$$
and
$$
\vert d\langle g_1\rangle\ldots d\langle g_n\rangle\vert_v\,=\,d_\Gamma(g_1,g_2)+\ldots+d_\Gamma(g_{n-1},g_n)+d_\Gamma(g_n,vg_1),
\eqno(4.8)
$$
respectively. The weight is invariant under left translation by $Z_v$. 
\item[c)] 
For $\lambda,N>1$ and an integer $k\geq 0$ let $\parallel-\parallel_{(\lambda,N,k)}$
be the largest norm on $CC_*({\mathbb C}\langle\Gamma\rangle)$ given on cyclic $n$-simplices by
$$
\parallel\omega_n\parallel_{(\lambda,N,k)}\,
=\,\frac{(1+n)^kN^{-n}}{c(n)!}\lambda^{\vert\omega_n\vert_v}.
\eqno(4.9)
$$
where $c(2n)=c(2n+1)=n$. 
\item[d)]
The linear map 
$$
\begin{array}{cccc}
\chi_v: & CC_*({\mathbb C}\langle\Gamma\rangle) & \longrightarrow & CC_*({\mathbb C}\langle\Gamma\rangle\rtimes U:{\mathbb C}\rtimes U)_{\{v\}} \\
 & & & \\
 & \langle g_0\rangle d\langle g_1\rangle\ldots d\langle g_n\rangle & \mapsto & u_v\langle g_0\rangle d(u_e\langle g_1\rangle)\ldots d(u_e\langle g_n\rangle) \\
 & & & \\
&  d\langle g_1\rangle\ldots d\langle g_n\rangle & \mapsto & d(u_v\langle g_1\rangle)d(u_e\langle g_2\rangle)\ldots d(u_e\langle g_n\rangle)\\
\end{array}
\eqno(4.10)
$$ 
is invariant under the $U$-action on the left and induces an isomorphism
$$
CC_*({\mathbb C}\langle\Gamma\rangle)^U\simeq CC_*({\mathbb C}\langle\Gamma\rangle)_U\,\overset{\simeq}{\longrightarrow}\,
CC_*({\mathbb C}\langle\Gamma\rangle\rtimes U:{\mathbb C}\rtimes U)_{\{v\}}.
\eqno(4.11)
$$ 
\item[e)]
The seminorms (4.9), restricted to the subcomplex\\ $CC_*({\mathbb C}\langle\Gamma\rangle)^U$ of $U$-invariants of $CC_*({\mathbb C}\langle\Gamma\rangle)$, give therefore rise, via (4.11), to seminorms on 
$CC_*({\mathbb C}\langle\Gamma\rangle\rtimes U:{\mathbb C}\rtimes U)_{\{v\}} $, denoted by the same letters.
\item[f)]
Let $CC_*({\mathbb C}\langle\Gamma\rangle\rtimes U:\,{\mathbb C}\rtimes U)_{(v,\lambda,N)}$ be the completion of \\
 $CC_*({\mathbb C}\langle\Gamma\rangle\rtimes U:\,{\mathbb C}\rtimes U)_{\{v\}}$ with respect to $\parallel-\parallel_{(\lambda,N,k)},\,k\geq 0$.\\ It is a Fr\'echet chain complex of $Z_v$-modules. 
\item[g)]
Introduce finally the ind-Fr\'echet complex of $Z_v$-modules 
$$
{\mathcal C}{\mathcal C}'_*({\mathbb C}\langle\Gamma\rangle\rtimes U:\,{\mathbb C}\rtimes U)_{\{v\}}\,=\,
"\underset{(\lambda,N)\to(1,\infty)}{\lim}"\,CC_*({\mathbb C}\langle\Gamma\rangle\rtimes U:\,{\mathbb C}\rtimes U)_{(v,\lambda,N)}.
\eqno(4.12)
$$
\end{itemize}
\end{deflem}

\begin{deflem} 
Let $\Gamma$ be a discrete isometry group of a $CAT(0)$-space as in Section 2 or let $\Gamma$ be word-hyperbolic. Let $v\in\Gamma$ be an element of finite order. 
The notations of 4.1. and 4.2 are understood.
\begin{itemize}
\item[a)] 
For real numbers $\lambda,N>1,$ and integers $k,l\geq 0$ let $\parallel-\parallel_{(\lambda,N,k,l)}$\\
be the largest norm on $C_*(\Gamma,v)\otimes CC_*({\mathbb C}\langle \Gamma\rangle)$
 given on standard generators  $\alpha_m\otimes\omega_n,\,\alpha_m\in\Delta_m(\Gamma),\,\omega_n\in\Delta^{cyc}_n(\Gamma)$
by
$$
\parallel\alpha_m\otimes\omega_n\parallel_{(\lambda,N,k,l)}\,=\,
\left(1+d_\Gamma(Supp(\alpha_m), Supp(\omega_n))+\vert\omega_n\vert_v\right)^l\cdot\parallel\omega_n\parallel_{(\lambda,N,k)}.
\eqno(4.13)
$$
\item[b)] The norms of a), restricted to the subcomplex $C_*(\Gamma,v)\otimes CC_*({\mathbb C}\langle \Gamma\rangle)^U$ yield, via (4.11), norms $\parallel-\parallel_{(\lambda,N,k,l)}$ on 
the complex $$C_*(\Gamma,v)\otimes CC_*({\mathbb C}\langle\Gamma\rangle\rtimes U:{\mathbb C}\rtimes U)_{\{v\}}.$$ 
\item[c)]
For $R\geq 0,\lambda,N>1$ let $CC_*^{tilt}(\Gamma,v)_{(R,\lambda,N)}$ be the completion of 
$$
C_*^{R}(\Gamma,v)^{alt}\otimes CC_*({\mathbb C}\langle \Gamma\rangle\rtimes U:\,{\mathbb C}\rtimes U)_{\{v\}}
$$
with respect to the seminorms $\parallel-\parallel_{(\lambda,N,k,l)},k,l\in{\mathbb N}$. This is a ${\mathbb Z}/2{\mathbb Z}$-graded complex of Fr\'echet-spaces.
\item[d)] Put
$$
{\mathcal C}{\mathcal C}_*^{tilt}(\Gamma,v)\,=\,"\underset{(R,\lambda,N)\to(\infty,1,\infty)}{\lim}"\,CC_*^{tilt}(\Gamma,v)_{(R,\lambda,N)}.
\eqno(4.14)
$$
This is a  ${\mathbb Z}/2{\mathbb Z}$-graded ind-complex of Fr\'echet-spaces.
\item[e)] The diagonal action of $Z_v$ on $C_*(\Gamma,v)\otimes CC_*({\mathbb C}\langle \Gamma\rangle)^U$ is isometric with respect to the norms $\parallel-\parallel_{(\lambda,N,k,l)}$ and descends to an isometric action on the Fr\'echet complexes $CC_*^{tilt}(\Gamma,v)_{(R,\lambda,N)}$ and on the ind-Fr\'echet complex  ${\mathcal C}{\mathcal C}_*^{tilt}(\Gamma,v)$.
\end{itemize}
\end{deflem}

\begin{proof}
It has to be shown that the differentials in the complexes are bounded with respect to the associated norms.
The differential of the Bar-complex is bounded because it is bounded degree-wise in the $\ell^1$-norm, and the complex $C_*^{R}(\Gamma,{\mathbb C})^{alt}$ equals zero in large degrees. 
One easily verifies that the Hochschild- and Connes-differentials of the Fr\'echet complexes $CC_*({\mathbb C}\langle\Gamma\rangle\rtimes U:\,{\mathbb C}\rtimes U)_{(\lambda,N)},\lambda,N>1,$ are bounded as well. 
\end{proof}

\begin{remark}
Already in the simplest case $v=e$ the complex ${\mathcal C}{\mathcal C}'_*({\mathbb C}\langle\Gamma\rangle)$ defined in 4.2 is rather different from the known analytic cyclic complexes as its topology is not related to any topology on the underlying algebra ${\mathbb C}\langle \Gamma\rangle$.
\end{remark}

\subsection{The basic retractions}

\begin{lemma}
The $Z_v$-equivariant chain map
$$
\begin{array}{ccc}
C_*(\Gamma,v)^{alt}\otimes CC_*({\mathbb C}\langle\Gamma\rangle
)^U & \overset{\epsilon\otimes id}{\longrightarrow} & CC_*({\mathbb C}\langle\Gamma\rangle\rtimes U:{\mathbb C}\rtimes U)_{\{v\}}
\end{array}
\eqno(4.15)
$$
is a chain homotopy equivalence of complexes of $Z_v$-modules.
\end{lemma}

\begin{proof}
We begin with the case of an isometry group of a $CAT(0)$-space.
The mapping cone of $\epsilon: C_*^{Bar}(X,v)^{alt}\to{\mathbb C}$ equals (up to a shift of degrees) the augmented Bar-complex 
 $\widetilde{C}_*^{Bar}(X,v)^{alt}=C_*^{Bar}(X,v)^{alt}\oplus{\mathbb C}[-1]$ with differential $\partial=\partial_{Bar}$ in strictly positive degrees and $\partial=\epsilon$ in degree zero. Extend the homotopy operators (2.21) , $y\in X^v$, to a contracting chain homotopy 
 $$
 h_{(y,v)}^{alt}: \widetilde{C}_*(\Gamma,v)^{alt}\,\to\, \widetilde{C}_{*+1}(\Gamma,v)^{alt},
 \eqno(4.16)
 $$
 of the augmented Bar-complex by putting $h_{(y,v)}^{alt}(1)=\pi_{(y,v)}([e])$ in degree -1.
 Define a strictly $Z_v$-equivariant linear map 
 $$
 \mu'_*:\,\widetilde{C}_*(\Gamma,v)^{alt}\otimes CC_*({\mathbb C}\langle\Gamma\rangle)\,
 \longrightarrow\,\widetilde{C}_{*+1}(\Gamma,v)^{alt}\otimes CC_*({\mathbb C}\langle\Gamma\rangle)
 \eqno(4.17)
 $$
 by
 $$
 \begin{array}{cc}
 \alpha_m\otimes\omega'_n  \mapsto 
  h^{alt}_{(\pi'_{X^v}(g_0x),v)}(\alpha_m)\otimes\omega'_n, &
 \alpha_m\otimes \omega''_n  \mapsto 
 h^{alt}_{(\pi'_{X^v}(g_1x),v)}(\alpha_m)\otimes\omega''_n\\ 
 \end{array}
 \eqno(4.18)
 $$
\\ 
for $\omega'_n=\langle g_0\rangle d\langle g_1\rangle\ldots d\langle g_n\rangle,\,\omega''_n= d\langle g_1\rangle\ldots d\langle g_n\rangle$ and $\alpha_m\in\Delta_m^R(\Gamma,v),\,m\geq -1$. 
\\
Let 
$$
\pi_U:CC_*({\mathbb C}\langle\Gamma\rangle)\to CC_*({\mathbb C}\langle\Gamma\rangle)^U
$$ 
be the canonical projection onto the $U$-fixed points given by averaging over the finite cyclic group $U$. 
The identification (4.11) allows to define the $Z_v$-equivariant linear map
$$
\mu=(id\otimes\chi_v)\circ (id\otimes\pi_U)\circ\mu'\circ (id\otimes\chi_v^{-1}):
\eqno(4.19)
$$
$$
\widetilde{C}_*(\Gamma,v)^{alt}\otimes CC_*({\mathbb C}\langle\Gamma\rangle\rtimes U:{\mathbb C}\rtimes U)_{\{v\}}\,\to\,
\widetilde{C}_{*+1}(\Gamma,v)^{alt}\otimes CC_*({\mathbb C}\langle\Gamma\rangle\rtimes U:{\mathbb C}\rtimes U)_{\{v\}}.
$$
Let
$$
\varphi\,=\,id-(\mu\circ\partial+\partial\circ\mu):
\eqno(4.20)
$$
$$
\widetilde{C}_*(\Gamma,v)^{alt}\otimes CC_*({\mathbb C}\langle\Gamma\rangle\rtimes U:{\mathbb C}\rtimes U)_{\{v\}}\,\to\,
\widetilde{C}_{*}(\Gamma,v)^{alt}\otimes CC_*({\mathbb C}\langle\Gamma\rangle\rtimes U:{\mathbb C}\rtimes U)_{\{v\}}
$$
This is a $Z_v$-equivariant chain map of the underlying ${\mathbb Z}/2{\mathbb Z}$-graded complexes which is $Z_v$-equivariantly chain homotopic to the identity
 and satisfies 
$$
\begin{array}{c}
\varphi\left(\widetilde{C}_*(\Gamma,v)^{alt}\otimes CC_*({\mathbb C}\langle\Gamma\rangle\rtimes U:{\mathbb C}\rtimes U)_{\{v\}}\right) \\
 \cap \\
\widetilde{C}_{*+1}(\Gamma,v)^{alt}\otimes CC_*({\mathbb C}\langle\Gamma\rangle\rtimes U:{\mathbb C}\rtimes U)_{\{v\}} \\
\end{array}
\eqno(4.21)
$$
for all $n\in{\mathbb N}$.\\
\\
The alternating Bar-complex yields a projective resolution of the constant $\Gamma$-module ${\mathbb C}$. As $\Gamma$ acts cocompactly on $X$ this module possesses also a projective resolution of finite length. As any two resolutions are chain homotopy equivalent it follows that there exists a $\Gamma$-equivariant linear map  $\eta':\widetilde{C}_*(\Gamma,{\mathbb C})^{alt}\to \widetilde{C}_{*+1}(\Gamma,{\mathbb C})^{alt}$ such that $\partial\circ\eta'+\eta'\circ\partial=id$ in sufficiently large degrees $*\geq N_0$. Put
$$
\eta=\pi_{X^v}\circ\eta'+h_{Bar}(\pi_{X^v},id):\widetilde{C}_*(\Gamma,v)^{alt}\to \widetilde{C}_{*+1}(\Gamma,v)^{alt}.
\eqno(4.22)
$$
Then
$$
h_0^{tilt}\,=\,(\eta\otimes id)\circ\varphi^{N_0+1}\,+\,\underset{k=0}{\overset{N_0}{\sum}}\,\mu\circ\varphi^k
\eqno(4.23)
$$
is a $Z_v$-equivariant contracting chain homotopy of 
$$
\widetilde{C}_*(\Gamma,v)^{alt}\otimes CC_*({\mathbb C}\langle\Gamma\rangle\rtimes U:{\mathbb C}\rtimes U)_{\{v\}}\,=\,Cone(\epsilon\otimes id)[1].
$$
In the case of a word-hyperbolic group the same argument applies with the following modifications. We set $h_{(y)}(1)=[y]$ in degree -1 and put
$$
\begin{array}{cccc}
 \mu': &\widetilde{C}_*(\Gamma,v)^{alt}\otimes CC_*({\mathbb C}\langle\Gamma\rangle) & \longrightarrow & 
\widetilde{C}_{*+1}(\Gamma,v)^{alt}\otimes CC_*({\mathbb C}\langle\Gamma\rangle) \\
 & & & \\
 & \alpha_m\otimes\omega'_n & \mapsto & \frac{1}{\vert S'_v(g_0)\vert}\underset{h\in S'_v(g_0)}{\sum}\,
  h^{alt}_{(h)}(\alpha_m)\otimes\omega'_n, \\
   & & & \\
  & \alpha_m\otimes \omega''_n & \mapsto & 
   \frac{1}{\vert S'_v(g_1)\vert}\underset{h'\in S'_v(g_1)}{\sum}\,
  h^{alt}_{(h')}(\alpha_m)\otimes\omega''_n, \\ 
 \end{array}
 \eqno(4.24)
$$
where 
$$
S'_v(g)\,=\,\{h\in Z_v,\,d(g,h)=d(g,Z_v)\},\,\,\,\forall g\in\Gamma.
\eqno(4.25)
$$
The centralizer $Z_v$ of an element of a word-hyperbolic group is word-hyperbolic itself \cite{BH} III.$\Gamma$, 4.7, 4.9. The alternating Rips-complexes $C_*^R(Z_v,{\mathbb C})^{alt}$ 
provide then finite resolutions of $Z_v$ by finitely generated projective $Z_v$-modules for $R>>0$.  We let $\eta=\eta':\widetilde{C}_*(\Gamma,v)^{alt}\to \widetilde{C}_*(\Gamma,v)^{alt}$ 
be any $Z_v$-equivariant linear map such that $\partial\circ\eta+\eta\circ\partial=id$ in sufficiently high degrees.
\end{proof}

\begin{lemma}
Let
$$
\left(\widetilde{C}_*^R(\Gamma,v)^{alt}\otimes CC_*({\mathbb C}\langle \Gamma\rangle),\parallel-\parallel_{(\lambda,N,k,l)}\right)\,=
$$
$$
=\,\left(C_*^R(\Gamma,v)^{alt}\otimes CC_*({\mathbb C}\langle \Gamma\rangle),\parallel-\parallel_{(\lambda,N,k,l)}\right)
\oplus \left(CC_*({\mathbb C}\langle \Gamma\rangle),\parallel-\parallel_{(\lambda,N,k)}\right)
$$
\begin{itemize}
\item[a)] In the $CAT(0)$ case the $Z_v$-equivariant linear map (4.17) satisfies
$$
\mu'(\widetilde{C}_*^R(\Gamma,v)^{alt}\otimes CC_*({\mathbb C}\langle \Gamma\rangle))\,\subset\,C_*^{R+2\,diam(D)+1}(\Gamma,v)^{alt}\otimes CC_*({\mathbb C}\langle \Gamma\rangle)
\eqno(4.26)
$$
and
$$
\parallel\mu'(\xi)\parallel_{(\lambda,N,k,l)}\,\leq\,C(\vert v\vert,R,diam(D),l)\parallel\xi\parallel_{(\lambda,N,k,l+1)}
\eqno(4.27)
$$
for all $\xi\in
\widetilde{C}_*^R(\Gamma,v)^{alt}\otimes CC_*({\mathbb C}\langle \Gamma\rangle)$.

\item[b)] In the $CAT(0)$ case the  $Z_v$-equivariant linear map (4.19) satisfies
$$
\mu(\widetilde{C}_*^R(\Gamma,v)^{alt}\otimes CC_*({\mathbb C}\langle\Gamma\rangle\rtimes U:{\mathbb C}\rtimes U)_{\{v\}})
$$
$$
\subset\,C_{*+1}^{R+2\,diam(D)+1}(\Gamma,v)^{alt}\otimes CC_*({\mathbb C}\langle\Gamma\rangle\rtimes U:{\mathbb C}\rtimes U)_{\{v\}}
\eqno(4.28)
$$
and
$$
\parallel\mu(\xi')\parallel_{(\lambda,N,k,l)}\,\leq\,C'(\vert v\vert,R,diam(D),l)\parallel\xi'\parallel_{(\lambda,N,k,l+1)},
\eqno(4.29)
$$
for all $\xi'\in
\widetilde{C}_*^R(\Gamma,v)^{alt}\otimes CC_*({\mathbb C}\langle\Gamma\rangle\rtimes U:{\mathbb C}\rtimes U)_{\{v\}}.$

\item[c)] In the $CAT(0)$ case there exists for each $R>0$ an $R'>>R$ such that the $Z_v$-equivariant linear map (4.22) satisfies
$$
(\eta\otimes id)(C_*^R(\Gamma,v)^{alt}\otimes CC_*({\mathbb C}\langle \Gamma\rangle))\,\subset\,C_*^{R'}(\Gamma,v)^{alt}\otimes CC_*({\mathbb C}\langle \Gamma\rangle)
\eqno(4.30)
$$
and
$$
\parallel(\eta\otimes id)(\xi)\parallel_{(\lambda,N,k,l)}\,\leq\,C(\eta,R)\parallel\xi\parallel_{(\lambda,N,k,l)}
\eqno(4.31)
$$
for all $\xi\in
\widetilde{C}_*^R(\Gamma,v)^{alt}\otimes CC_*({\mathbb C}\langle \Gamma\rangle).$
\end{itemize}
Similar estimates hold in the case of a word-hyperbolic group with universal constants depending now on $[v],R,l$ and $\delta$.
\end{lemma}

\begin{proof}
The norms $\parallel -\parallel_{(\lambda,N,k,l)},\lambda,N>1,k,l\in{\mathbb N},$ are weighted $\ell^1$-norms on the complex $Cone(\epsilon\otimes id)$ spanned by the bisimplicial set $(\Delta_\bullet(\Gamma)\cup\{1\})\times\Delta^{cyc}_\bullet(\Gamma).$  
It suffices therefore to majorize the norm on bi-simplices. \\
\\
a) We find for $\alpha_m=[h_0,\ldots,h_m]\in\Delta_m^R(\Gamma,v),\,\omega'_n=\langle g_0\rangle d\langle g_1\rangle\ldots d\langle g_n\rangle\in\Delta_n^{cyc}(\Gamma)'$
$$
\mu'(\alpha_m\otimes\omega'_n)\,=\,h^{alt}_{(\pi'_{X^v}(g_0x),v)}(\alpha_m)\otimes\omega'_n\,=\,
\left(\pi_{alt}\circ\underset{k=0}{\overset{\infty}{\sum}}\,h^{Bar}(\pi_{(\pi'_{X^v}(g_0x),v,k)},\pi_{(\pi'_{X^v}(g_0x),v,k+1)})(\alpha_m)\right)\otimes\omega'_n
$$
where
$$
h^{Bar}(\pi_{(\pi'_{X^v}(g_0x),v,k)},\pi_{(\pi'_{X^v}(g_0x),v,k+1)})(\alpha_m)\,=\,\underset{i=0}{\overset{m}{\sum}}\,(-1)^i\alpha_m^{(k,i)}
\eqno(4.32)
$$
with 
$$
\alpha_m^{(k,i)}\,=\,[\pi_{(\pi'_{X^v}(g_0x),v,k)}(h_0),\ldots,\pi_{(\pi'_{X^v}(g_0x),v,k)}(h_i),\pi_{(\pi'_{X^v}(g_0x),v,k+1)}(h_i),
\ldots,\pi_{(\pi'_{X^v}(g_0x),v,k+1)}(h_m)].
$$
Note that $d(g,g')\leq d(h,h')+\vert k-k'\vert+2\,diam(D)$ for all $g\in Supp(\pi_{(\pi'_{X^v}(g_0x),v,k)}(h))$ and $g'\in Supp(\pi_{(\pi'_{X^v}(g_0x),v,k')}(h'))$ by (2.17)  so that 
$$
h^{alt}_{(\pi'_{X^v}(g_0x),v)}(\alpha_m)\in C_*^{R+2\,diam(D)+1}(\Gamma,v).
\eqno(4.33)
$$
We estimate now the distance of $\alpha_m^{(k,i)}\otimes\omega'_n$ from the diagonal.\\ For 
$g\in Supp(\pi_{(\pi'_{X^v}(g_0x),v,k)}(h_0))$ we find on the one hand
$$
d(gx,\pi'_{X^v}(g_0x))\,\leq\,d(\pi'_{X^v}\circ\pi'_{B(\pi'_{X^v}(g_0x),k)}(h_0x),\pi'_{X^v}(g_0x))+diam(D)
$$
$$
\leq d(h_0x,g_0x)+diam(D)
$$
and on the other hand $d(\pi'_{X^v}(g_0x),g_0x)\,\leq\,d(h_0x,g_0x)+diam(D)$\\
as $h_0D\cap X^v\neq\emptyset$, so that 
$$
d(Supp(\alpha_m^{(k,i)}),Supp(\omega'_n))\,\leq\,d(gx,g_0x)\,\leq\,2d(h_0x,g_0x)+2\,diam(D).
\eqno(4.34)
$$
Concerning $\ell^1$-norms we note that 
$$
\parallel h_{Bar}(\pi_{(\pi'_{X^v}(g_0x),v,k)},\pi_{(\pi'_{X^v}(g_0x),v,k+1)})(\alpha_m)\parallel_{\ell^1(\Delta_\bullet(\Gamma))}\,\leq\,m+1
\eqno(4.35)
$$
As $h_{Bar}(\pi_{(\pi'_{X^v}(g_0x),v,k)},\pi_{(\pi'_{X^v}(g_0x),v,k+1)})(\alpha_m)$ vanishes as soon as 
$$
Supp(\alpha_m)\cdot x\,\subset\,B(\pi'_{X^v}(g_0x),k)
$$
one only gets contributions to the sum (4.32) for
$$
k\,\leq\,d(\pi'_{X^v}(g_0x),h_0x)+R\,\leq\,d(\pi'_{X^v}(g_0x),\pi'_{X^v}(h_0x))+d(\pi'_{X^v}(h_0x),h_0x)+R
$$
$$
\leq d(g_0x,h_0x)+diam(D)+R.
\eqno(4.36)
$$
Finally 
$$
d(h_0x,g_0x)\,\leq d(Supp(\alpha_m),Supp(\omega'_n))+diam(Supp(\alpha_m))+diam(Supp(\omega'_n))
$$
$$
\leq\,d(Supp(\alpha_m),Supp(\omega'_n))+R+\vert\omega'_n\vert_v.
\eqno(4.37)
$$
Altogether this shows that
$$
\mu':\,C_*^R(\Gamma,v)^{alt}\otimes CC_*({\mathbb C}\langle \Gamma\rangle)\to C_{*+1}^{R+2\,diam(D)+1}(\Gamma,v)^{alt}\otimes CC_*({\mathbb C}\langle \Gamma\rangle)
$$
satisfies
$$
\parallel\mu'(\alpha_m\otimes\omega'_n)\parallel_{(\lambda,N,k,l)}\,\leq\,C(R,diam(D),l,m)\cdot\parallel\alpha_m\otimes\omega'_n\parallel_{(\lambda,N,k,l+1)}
$$
as was to be shown. The proofs of the remaining cases are similar.
For hyperbolic groups one uses the fact that $Z_v$ is a quasi-convex subset of ${\mathcal G}(\Gamma,S)$ \cite{BH},III,$\Gamma$,3.9.\\
\\
b)
By construction the linear maps $id\otimes\chi_v^{-1} $ and $id\otimes\chi_v$ are isometric with respect to the given seminorms. Concerning the averaging operator $id\otimes\pi_U$ we observe for 
the generator $v\in U$ 
$$
\parallel (id\otimes v)(\alpha_m\otimes\omega'_n)\parallel_{(\lambda,N,k,l)}\,=\,
\left(1+d(Supp(\alpha_m),Supp(v\omega'_n))+\vert v\omega'_n\vert_v\right)^l\cdot\parallel v\omega'_n\parallel_{(\lambda,N,k)}
$$
$$
\leq\,\left(1+d(Supp(\alpha_m),Supp(\omega'_n))+2\vert\omega'_n\vert_v\right)^l\cdot\parallel\omega'_n\parallel_{(\lambda,N,k)}
\,\leq\,
2^l\cdot\parallel \alpha_m\otimes\omega'_n\parallel_{(\lambda,N,k,l)}
$$
which follows from $(g,vg)\leq\vert\omega'_n\vert_v,\,\forall g\in Supp(\omega'_n),$ and the $Z_v$-invariance of the norms $\parallel-\parallel_{(\lambda,N,k)}$. 
Together with part a) this shows our claim.\\
\\
c) The $\Gamma$-equivariance of $\eta':\,C_*(\Gamma,{\mathbb C})^{alt}\to C_{*+1}(\Gamma,{\mathbb C})^{alt}$ and the fact that 
$C_*^R(\Gamma,{\mathbb C})^{alt}$ is a finitely generated projective $\Gamma$-module imply that 
$$
\eta(C_*^R(\Gamma,{\mathbb C})^{alt})\subset C_{*+1}^{R'}(\Gamma,{\mathbb C})^{alt}
$$ 
for $R'>>R$ large enough. It is easily seen 
that $\eta\otimes id$ is bounded with respect to each of the seminorms $\parallel-\parallel_{(\lambda,N,k,l)},\,\lambda,N>1,k,l\in{\mathbb N}$.
\end{proof}

\begin{lemma}
The $Z_v$-equivariant chain map
$$
\begin{array}{ccc}
C_*(\Gamma,v)^{alt}\otimes CC_*({\mathbb C}\langle\Gamma\rangle\rtimes U:{\mathbb C}\rtimes U)_{\{v\}} & \overset{id\otimes CC(\epsilon)}{\longrightarrow} & 
C_*(\Gamma,v)^{alt}\otimes CC_*({\mathbb C}\langle\{*\}\rangle\rtimes U:{\mathbb C}\rtimes U)_{\{v\}} \\
\end{array}
\eqno(4.38)
$$
is a chain homotopy equivalence in the category of complexes of $Z_v$-modules.\\ Here $\{*\}$ is the one point set equipped with the trivial $Z_v$-action.
\end{lemma}

\begin{proof}
Recall that the cone of a chain map $f:C\to C'$ of complexes is given by the complex
$Cone(f)=C[1]\oplus C'$ with differential 
$$
\partial_{Cone}\,=\,\left(
\begin{matrix}
-\partial_C & 0 \\
f & \partial_{C'} \\
\end{matrix}
\right)
\eqno(4.39)
$$
and that $f$ is a chain-homotopy equivalence iff $Cone(f)$ is contractible.
We consider the chain map 
$$
id\otimes CC(\epsilon):\,C_*(\Gamma,v)\otimes CC_*({\mathbb C}\langle\Gamma\rangle)\,\longrightarrow\,
C_*(\Gamma,v)\otimes CC_*({\mathbb C}\langle\{*\}\rangle)
$$
of complexes of $Z_v$-modules and define a $Z_v$-equivariant linear operator
$$
\nu=\left(
\begin{matrix}
h' & j' \\
0 & 0 \\
\end{matrix}
\right):\,Cone(id\otimes CC(\epsilon))_*\longrightarrow\,Cone(id\otimes CC(\epsilon))_{*+1}
\eqno(4.40)
$$
by putting (in the notations of 3.4)
$$
h'(\alpha_m\otimes\omega_n)\,=\,(-1)^{m-1}\pi_{alt}(\alpha_m)\otimes\frac{1}{\vert U\vert}\cdot\underset{u\in U}{\sum}\,h_{uh_0}^{cyc}(\omega_n)
\eqno(4.41)
$$
and
$$
j'(\alpha_m\otimes\widetilde{\omega})\,=\,\pi_{alt}(\alpha_m)\otimes\frac{1}{\vert U\vert}\cdot\underset{u\in U}{\sum}\,(\pi_{uh_0}\circ i)_*(\widetilde{\omega})
\eqno(4.42)
$$
for $\alpha_m=[h_0,\ldots,h_m]\in\Delta_m(\Gamma),\,\omega_n\in\Delta_n^{cyc}({\mathbb C}\langle\Gamma\rangle)$, and $\widetilde{\omega}\in\Delta_*^{cyc}({\mathbb C}\langle\{*\}\rangle)$.
Here $j:\{*\}\to \Gamma$ sends $*$ to $e$. The $Z_v$-equivariant chain map
$$
\psi=id-(\partial\circ\nu+\nu\circ\partial):Cone(id\otimes CC(\epsilon))\longrightarrow\,Cone(id\otimes CC(\epsilon))
\eqno(4.43)
$$
satisfies then
$$
\psi\left(C_*(\Gamma,v)\otimes Cone(CC(\epsilon))_*\right)\,\subset\,C_{*-1}(\Gamma,v)\otimes Cone(CC(\epsilon))_*,
\eqno(4.44)
$$
so that 
$$
h_1^{tilt}\,=\,\underset{n=0}{\overset{\infty}{\sum}}\,\nu\circ\psi^n:\,Cone(id\otimes CC(\epsilon))_*\longrightarrow\,Cone(id\otimes CC(\epsilon))_{*+1}
\eqno(4.45)
$$
is a well defined $Z_v$-equivariant contracting chain homotopy of $Cone(id\otimes CC(\epsilon))_*$.
By construction $h_1^{tilt}$ preserves the subcomplex given by the cone of 
$$
id\otimes CC(\epsilon):\,C_*(\Gamma,v)^{alt}\otimes CC_*({\mathbb C}\langle\Gamma\rangle)^U\,\longrightarrow\,
C_*(\Gamma,v)^{alt}\otimes CC_*({\mathbb C}\langle\{*\}\rangle)
$$
so that it provides, via (4.11), a $Z_v$-equivariant contracting homotopy of the cone of
$$
C_*(\Gamma,v)^{alt}\otimes CC_*({\mathbb C}\langle\Gamma\rangle\rtimes U:{\mathbb C}\rtimes U)_{\{v\}}\,\overset{id\otimes CC(\epsilon)}{\longrightarrow}\,
C_*(\Gamma,v)^{alt}\otimes CC_*({\mathbb C}\langle\{*\}\rangle\rtimes U:{\mathbb C}\rtimes U)_{\{v\}}
$$
 as claimed.
\end{proof}

\begin{lemma}
The chain map
$$
\widehat{CC}_*({\mathbb C}\langle\{*\}\rangle\rtimes U:{\mathbb C}\rtimes U)_{\{v\}}\,\longrightarrow\,{\mathbb C}
\eqno(4.46)
$$
which sends $u_v\langle *\rangle$ to 1 and vanishes on $\Omega^{\geq 1}({\mathbb C}\langle\{*\}\rangle\rtimes U:{\mathbb C}\rtimes U)_\natural$ is a chain homotopy equivalence.
Here $\widehat{CC}_*=\left(\underset{n=0}{\overset{\infty}{\prod}}\,\Omega^n_\natural,b+B\right)$ is the periodic cyclic bicomplex.
\end{lemma}

A chain homotopy inverse is given by the linear map sending $1\in{\mathbb C}$ to $u_v\cdot ch(\langle *\rangle)\in \widehat{CC}_0({\mathbb C}\langle\{*\}\rangle\rtimes U:{\mathbb C}\rtimes U)_{\{v\}},$ where $ch$ denotes the Chern-character \cite{Lo} 
of the idempotent $\langle *\rangle\in{\mathbb C}\langle\{*\}\rangle$.\\
\\
We let
$$
\epsilon':\, CC_*({\mathbb C}\langle\Gamma\rangle\rtimes U:{\mathbb C}\rtimes U)_{\{v\}}\,\longrightarrow\,
CC_*({\mathbb C}\langle\{*\}\rangle\rtimes U:{\mathbb C}\rtimes U)_{\{v\}}\,\longrightarrow\,{\mathbb C}
\eqno(4.47)
$$
be the composition of the $Z_v$-equivariant chain maps (4.38) and (4.46).

\begin{lemma}
In the notations of 4.7 we have the estimates
$$
\parallel h'(\alpha_m\otimes\omega_n)\parallel_{(\lambda,\lambda^{2\,diam(D)+1}N,k,l)}\,\leq\,C(\lambda,R,diam(D),k,l)\parallel\alpha_m\otimes\omega_n\parallel_{(\lambda,N,k+l+2,l+1)}
\eqno(4.48)
$$
and
$$
\parallel j'(\alpha_m\otimes\widetilde{\omega}_n)\parallel_{(\lambda,N,k,l)}\,\leq\,C(\lambda,diam(D),l)\parallel\widetilde{\omega}_n\parallel_{(\lambda,N,k)}
\eqno(4.49)
$$
for suitable constants $C(\lambda,R,diam(D),k,l)$ and $C(\lambda,diam(D),l)$ and all\\ $\alpha_m\in\Delta_m^R(\Gamma,v),\,\omega_n\in\Delta_n^{cyc}({\mathbb C}\langle\Gamma\rangle),$ and $\widetilde{\omega}_n\in\Delta_n^{cyc}({\mathbb C}\langle *\rangle).$
\end{lemma}

\begin{proof}
Suppose at first that $\Gamma$ acts on a $CAT(0)$-space.\\
Let $\alpha_m=[h_0,\ldots,h_m],\omega_n=\langle g_0\rangle d\langle g_1\rangle\ldots d\langle g_n\rangle$ and let $k\in{\mathbb N}$. Then we find, according to (4.41) and (3.26)
$$
h'(\alpha_m\otimes\omega_n)\,=\,(-1)^{m-1}\pi_{alt}(\alpha_m)\otimes\frac{1}{\vert U\vert}\,\underset{u\in U}{\sum}\,\underset{d=0}{\overset{\infty}{\sum}}\,
(h+H)(\pi_{(uh_0x,e,d)},\pi_{(uh_0x,e,d+1)})(\omega_n)
$$
where $(h+H)(\pi_{(uh_0x,e,d)},\pi_{(uh_0x,e,d+1)})(\omega'_n)=H(\pi_{(uh_0x,e,d)},\pi_{(uh_0x,e,d+1)})(\omega'_n)$ is a linear combination of cyclic simplices 
$$
\omega'_n=d\langle g'_i\rangle\ldots d\langle g'_n\rangle d\langle g'_0\rangle\ldots d\langle g'_{i-1}\rangle
$$
with
$$
g'_j\,\in\,Supp(\pi_{(uh_0x,e,d)}(g_j))\cup Supp(\pi_{(uh_0x,e,d+1)}(g_j)),\,0\leq j\leq n.
$$
Thus $d_\Gamma(g'_j,g'_{j'})\,\leq\,d_\Gamma(g_j,g_{j'})+2\,diam(D)+1$ and
$$
\vert\omega'_n\vert_v\,\leq\,d(g_i,g_{i+1})+\ldots+d(g_n,g_0)+\ldots+d(g_{i-1},g_i)+d(g'_i,vg'_i)+(n+1)(2\,diam(D)+1)
$$
The convexity of the distance function and the fact that $uh_0D\cap X^v\neq\emptyset$ imply
$$
d(g'_i,vg'_i)\,\leq\,d(\pi'_{B(uh_0x,d')}(g_ix),v\pi'_{B(uh_0x,d')}(g_ix))+2\,diam(D)\,\leq\,d(g_i,vg_i)+d(uh_0x,vuh_0x)+2\,diam(D)
$$
so that
$$
\vert\omega'_n\vert_v\,\leq\,\vert\omega_n\vert_v+d(g_0,vg_0)+d(g_i,vg_i)+(n+3)(2\,diam(D)+1)\,\leq\,
3\vert\omega_n\vert_v+(n+3)(2\,diam(D)+1)
\eqno(4.50)
$$
The vertices of $\omega'_n$  lie at distance at most $diam(D)$ from a geodesic segment joining an element of $u\cdot Supp(\alpha_m)$ and an element 
of $Supp(\omega_n)$. Consequently
$$
d(Supp(\alpha_m),Supp(\omega'_n))\,\leq\,d(Supp(\alpha_m),Supp(\omega_n))+R+3\,diam(D)
\eqno(4.51)
$$
for all $h_i\in Supp(\alpha_m)$ and $u\in U$. For the $\ell^1$-norms one finds
$$
\parallel H(\pi_{(uh_0x,d)},\pi_{(uh_0x,d+1)})(\omega_n)\parallel_{\ell^1(\Delta_\bullet^{cyc}(\Gamma)}\,\leq\,n(n+1).
\eqno(4.52)
$$
Moreover $H(\pi_{(uh_0x,d)},\pi_{(uh_0x,d+1)})(\omega_n)=0$ unless 
$$
d< d(Supp(\alpha_m),Supp(\omega_n))+\vert\omega_n\vert_v+2\,diam(D).
\eqno(4.53)
$$
as $\pi_{(uh_0x,d)}(\omega_n)=\pi_{(uh_0x,d+1)}(\omega_n)=\omega_n$ otherwise. Altogether we obtain
$$
\parallel h'(\alpha_m\otimes\omega_n)\parallel_{(\lambda,\lambda^{2\,diam(D)+1}N,k,l)}\,\leq\,C(\lambda,R,diam(D),k,l)\parallel\alpha_m\otimes\omega_n\parallel_{(\lambda,N,k+l+2,l+1)}
$$
as claimed. For cyclic simplices $\omega_n\in\Delta_n^{cyc}(\Gamma)''$ the reasoning is similar. The second estimate (4.49) is trivial. The argument is almost verbatim the same 
for hyperbolic groups with constants depending on $(\lambda,R,\delta,\ell_S(v),k,l)$ and $(\lambda,\ell_S(v),l)$, respectively.
\end{proof}

\begin{theorem}
Let $\Gamma$ be a discrete group acting properly, isometrically and cocompactly on a $CAT(0)$ or suppose that $\Gamma$ is word-hyperbolic. Let $v\in\Gamma$ be an element of finite order and let $U\subset\Gamma$ be the subgroup generated by $v$. Then the chain maps

$$
\begin{array}{ccccc}
CC_*({\mathbb C}\langle\Gamma\rangle\rtimes U:{\mathbb C}\rtimes U)_{\{v\}}  & \overset{\epsilon\otimes id}{\longleftarrow} & C_*(\Gamma,v)^{alt}\otimes 
CC_*({\mathbb C}\langle\Gamma\rangle\rtimes U:{\mathbb C}\rtimes U)_{\{v\}}  & \overset{id\otimes \epsilon'}{\longrightarrow} & C_*(\Gamma,v)^{alt} 
\end{array}
$$
\\
induce bounded chain-homotopy equivalences

$$
\begin{array}{ccccc}
{\mathcal C}{\mathcal C}'_*({\mathbb C}\langle\Gamma\rangle\rtimes U:{\mathbb C}\rtimes U)_{\{v\}} & \overset{\sim}{\longleftarrow} & {\mathcal C}{\mathcal C}_*^{tilt}(\Gamma,v) & \overset{\sim}{\longrightarrow} &  {\mathcal C}^{alt}_*(\Gamma,v)
\end{array}
\eqno(4.54)
$$
\\
in the ind-category of Fr\'echet-complexes of $\Gamma$-modules.
\end{theorem}

\begin{proof}
Lemma 4.5 and 4.6 show that the linear endomorphism (4.23) of
$$
Cone\left(\epsilon\otimes id:C_*(\Gamma,v)^{alt}\otimes CC_*({\mathbb C}\langle\Gamma\rangle\rtimes U:\,{\mathbb C}\rtimes U)_{\{v\}}\,\longrightarrow\,
CC_*({\mathbb C}\langle\Gamma\rangle\rtimes U:\,{\mathbb C}\rtimes U)_{\{v\}}\right)
$$
gives rise to a $Z_v$-equivariant contracting chain homotopy of the  ind-Fr\'echet-complex 
$$
Cone\left(\epsilon\otimes id:\,{\mathcal C}{\mathcal C}_*^{tilt}(\Gamma,v)\,\longrightarrow\,{\mathcal C}{\mathcal C}'_*({\mathbb C}\langle\Gamma\rangle\rtimes U:\,{\mathbb C}\rtimes U)_{\{v\}}\right).
$$
This shows that the morphism
$$
\epsilon\otimes id:\,{\mathcal C}{\mathcal C}_*^{tilt}(\Gamma,v)\,\longrightarrow\,{\mathcal C}{\mathcal C}'_*({\mathbb C}\langle\Gamma\rangle\rtimes U:\,{\mathbb C}\rtimes U)_{\{v\}}
$$
is a chain homotopy equivalence in the category of ind-Fr\'echet-complexes of $Z_v$-modules.
Similarly lemma 4.7 and 4.9 show that the linear endomorphism (4.45) of the cone of 
$$
\begin{array}{c}
C_*(\Gamma,v)^{alt}\otimes CC_*({\mathbb C}\langle\Gamma\rangle\rtimes U:\,{\mathbb C}\rtimes U)_{\{v\}}\\
 \\
id\otimes CC(\epsilon)\,\downarrow \\
\\
C_*(\Gamma,v)^{alt}\otimes CC_*({\mathbb C}\langle *\rangle\rtimes U:\,{\mathbb C}\rtimes U)_{\{v\}}\\
\end{array}
$$
gives rise to a $Z_v$-equivariant contracting chain homotopy of the  ind-Fr\'echet-complex 
$$
Cone\left(id\otimes CC(\epsilon):\,{\mathcal C}{\mathcal C}_*^{tilt}(\Gamma,v)\,\longrightarrow\,{\mathcal C}_*^{alt}(\Gamma,v)\otimes_\pi
{\mathcal C}{\mathcal C}'_*({\mathbb C}\langle *\rangle\rtimes U:\,{\mathbb C}\rtimes U)_{\{v\}}\right).
$$
Together with 4.8 this shows that the composition 
$$
id\otimes\epsilon':{\mathcal C}{\mathcal C}_*^{tilt}(\Gamma,v)\,\longrightarrow\,{\mathcal C}_*^{alt}(\Gamma,v)\otimes_\pi
{\mathcal C}{\mathcal C}'_*({\mathbb C}\langle *\rangle\rtimes U:\,{\mathbb C}\rtimes U)_{\{v\}}\,\longrightarrow\,{\mathcal C}_*^{alt}(\Gamma,v)
$$
is a homotopy equivalence in the category of ind-Fr\'echet-complexes of $Z_v$-modules.
\end{proof}

\section{Analytic and local cyclic homology}

Before we come to our main results let us recall the definition of analytic and local cyclic cohomology of Banach algebras.
\begin{definition} \cite{Pu2}
Let $A$ be a complex Banach algebra. For $N>1$ and $k\in{\mathbb N}$ let $\parallel-\parallel_{(N,k)}$ be the largest norm on the cyclic bicomplex $CC_*(A)$ satisfying
$$
\parallel a_0da_1\ldots da_n\parallel_{(N,k)}\,\leq\,\frac{(1+n)^k}{c(n)!}N^{-n}\parallel a_0\parallel_A\cdot\ldots\cdot\parallel a_n\parallel_A
\eqno(5.1)
$$
and
$$
\parallel da_1\ldots da_n\parallel_{(N,k)}\,\leq\,\frac{(1+n)^k}{c(n)!}N^{-n}\parallel a_1\parallel_A\cdot\ldots\cdot\parallel a_n\parallel_A
\eqno(5.2)
$$
for all $a_0,a_1,\ldots,a_n\in A$ and all $n\in{\mathbb N}$, where $c(2n)=c(2n+1)=n$. The completion of the cyclic bicomplex $CC_*(A)$ with respect to the seminorms 
$\parallel-\parallel_{(N,k)},\,k\in{\mathbb N}$ is a Fr\'echet-chain complex denoted by $CC_*(A)_N$. The formal inductive limit of these complexes is denoted by
$$
CC^\omega_*(A)="\underset{N\to\infty}{\lim}"\,CC_*(A)_{N}
\eqno(5.3)
$$
\end{definition}
It is a natural ind-Fr\'echet-complex attached to $A$.\\

\begin{definition}
Let $A$ be a separable (not necessarily unital) Banach algebra. Denote by $Alg\to A$ the essentially small, filtrant category with compact homomorphisms $\varphi_{A'}:A'\to A$ of separable Banach algebras as objects and morphisms\\ $(A',\varphi_{A'})\to(A'',\varphi_{A''})$ given by homomorphisms $f:A'\to A''$ of Banach algebras compatible with the structure maps, i.e. satisfying $\varphi_{A''}\circ f=\varphi_{A'}$.\\
\\
The {\bf analytic cyclic bicomplex} of the Banach algebra $A$ is the ind-complex
$$
\begin{array}{ccccc}
{\mathcal C}{\mathcal C}_*(A) & = & \underset{Alg\to A}{\lim}CC_*^\omega(A') & = &  "\underset{Alg\to A}{\lim}\,\underset{N\to\infty}{\lim}CC_*(A')_N "\\
\end{array}
\eqno(5.4)
$$
Strictly speaking the formal inductive limit has to be taken over a small, cofinal, filtrant subcategory  of $Alg\to A$. It is, up to canonical isomorphism, independent of the choice of this subcategory.
The bivariant {\bf analytic cyclic cohomology} of a pair $(A,B)$ of separable complex Banach algebras is defined as
$$
HC^{an}_*(A,B)\,=\,Mor_{Ho(ind({\mathcal F}{\mathcal C}))}({\mathcal C}{\mathcal C}_*(A) ,{\mathcal C}{\mathcal C}_*(B) ),
\eqno(5.5)
$$
the vector space of chain-homotopy classes of morphisms in the category of ind-Fr\'echet-complexes between the analytic cyclic bicomplexes of $A$ and $B$, respectively.
\end{definition}

\begin{definition}
Let $Ho(ind({\mathcal F}{\mathcal C}))$ be the chain-homotopy category of ${\mathbb Z}/2{\mathbb Z}$-graded ind-Fr\'echet complexes. It is triangulated by declaring a triangle distinguished iff it is isomorphic to a 
triangle ${\mathcal C}'_*\,\overset{f}{\to}\,{\mathcal C}_*\,\to\,Cone(f)_*\,\to\,{\mathcal C}'_*[1]$. An ind-complex is constant if it is labeled by the filtrant category with a single morphism. An ind-complex ${\mathcal C}_*$ is 
{\bf weakly contractible} if every morphism from a constant ind-complex to it is null-homotopic. The weakly contractible ind-Fr\'echet-complexes form a null-system \cite{KS} in the triangulated chain homotopy category $Ho(ind({\mathcal F}{\mathcal C}))$. Its localization with respect to this null-system is called the {\bf derived ind-category} $ind({\mathcal D})$. The bivariant {\bf local cyclic cohomology} of the pair $(A,B)$ of complex, separable Banach algebras is defined as
$$
HC^{loc}_*(A,B)\,=\,Mor_{ind({\mathcal D})}({\mathcal C}{\mathcal C}_*(A) ,{\mathcal C}{\mathcal C}_*(B) ),
\eqno(5.6)
$$
the vector space of morphisms in the derived ind-category between the analytic cyclic bicomplexes of $A$ and $B$, respectively.
\end{definition}

Let $\Gamma$ be a finitely generated discrete group with word-length function $\ell$. For $\lambda>1$ we let $\ell^1_\lambda(\Gamma,\ell)$ be the completion of the group ring ${\mathbb C}\rtimes\Gamma$ with respect to the largest norm $\parallel-\parallel_\lambda$ satisfying $\parallel g\parallel_\lambda\leq\lambda^{\ell(g)}$ for all $g\in\Gamma$. It is a Banach algebra. The identity map induces a homomorphism $\ell^1_\lambda(\Gamma,\ell)\to\ell^1_{\lambda'}(\Gamma,\ell)$ for $\lambda<\lambda'$ and the topological direct limit of this family equals $\underset{\lambda\to 1}{\lim}\,\ell^1_\lambda(\Gamma,\ell)\,=\,\ell^1(\Gamma).$\\
\\
Let $v\in\Gamma$ and denote by $[v]\subset\Gamma$ its conjugacy class. The cyclic bicomplex of the group ring ${\mathbb C}\Gamma$ decomposes into a direct sum  
$$
CC_*({\mathbb C}\Gamma)\,=\,\underset{[v]}{\bigoplus}\,CC_*({\mathbb C}\Gamma)_{[v]},
\eqno(5.7)
$$
labeled by the conjugacy classes of $\Gamma$, of the subcomplexes $CC_*({\mathbb C}\Gamma)_{[v]}$ spanned by differential forms 
$g_0dg_1\ldots dg_n,\,g_0\ldots g_n\in[v]$ and $dh_1\ldots dh_n,\,h_1\ldots h_n\in[v],\,n\in{\mathbb N}$. This is called the {\bf homogeneous decomposition} of the cyclic bicomplex of a group ring. 
The contribution of the conjugacy classes of the unit, the torsion elements and the elements of infinite order are called the {\bf homogeneous, elliptic} and {\bf hyperbolic} part of the cyclic bicomplex, respectively. 
The homogeneous decomposition  gives rise to an isometric isomorphism of normed vector spaces
$$
\left(CC_*({\mathbb C}\Gamma),\parallel-\parallel_{(\lambda,N,k)}\right)\,\simeq\,\underset{[v]}{\bigoplus}\,\left(CC_*({\mathbb C}\Gamma)_{[v]},\parallel-\parallel_{(\lambda,N,k)}\right)
\eqno(5.8)
$$
and an isometric isomorphism
$$
\left(CC_*(\ell^1_\lambda(\Gamma))_N,\parallel-\parallel_{(N,k)}\right)\,\simeq\,\underset{[v]}{\widehat{\bigoplus}}\,\left((CC_*(\ell^1_\lambda(\Gamma))_{[v]})_{N},\parallel-\parallel_{(N,k)}\right)
\eqno(5.9)
$$
of Banach-spaces for all $\lambda,N>1$ and $k\in{\mathbb N}$. Consequently the ind-Fr\'echet complexes $CC^\omega_*(\ell^1_\lambda(\Gamma))$ decompose into a topological direct sum
$$
CC^\omega_*(\ell^1_\lambda(\Gamma))\,\simeq\,\underset{[v]}{\widehat{\bigoplus}}\,CC^\omega_*(\ell^1_\lambda(\Gamma))_{[v]},
\eqno(5.10)
$$
labeled by the conjugacy classes of $\Gamma$.

\subsection{The homogeneous and the elliptic part}
\begin{theorem}
Let $\Gamma$ be a discrete group acting properly, isometrically and cocompactly on a $CAT(0)$-space or suppose that $\Gamma$ is word-hyperbolic. Let $v\in\Gamma$ be an element of finite order and let $[v]\subset\Gamma$ be its conjugacy class. Then there exists a chain-homotopy equivalence 
$$
\underset{\lambda\to 1}{\lim}\,CC_*^\omega(\ell^1_\lambda(\Gamma))_{[v]}\,\overset{\sim}{\longrightarrow}\,H_*(\Gamma, {\mathbb C}[v])
\eqno(5.11)
$$
in the category of ind-Fr\'echet-complexes. Here the right hand side equals the constant, finite dimensional ind-complex with vanishing differentials given by the direct sum of the homology 
of even, respectively odd, degrees of $\Gamma$ with coefficient in the $\Gamma$-module ${\mathbb C}[v]$, equipped with the adjoint action.
\end{theorem}
.\\
\begin{proof}\\
{\bf Step 1:}\\
For a Fr\'echet $\Gamma$-module $M$ let $M_\Gamma$ be the associated space of $\Gamma$-coinvariants, i.e. the quotient of $M$ by the closure 
of the linear subspace spanned by the elements $g\xi-\xi,\,g\in\Gamma,\xi\in M$. The quotient norms on $M_\Gamma$ associated to a family of seminorms defining the Fr\'echet-structure on $M$ 
turn $M_\Gamma$ into a Fr\'echet space.\\
{\bf Step 2:}\\
Theorem 4.10 yields, after passing to coinvariants, chain homotopy equivalences
$$
\begin{array}{ccccc}
\left({\mathcal C}{\mathcal C}'_*({\mathbb C}\langle\Gamma\rangle\rtimes U:{\mathbb C}\rtimes U)_{\{v\}}\right)_{Z_v} & \overset{\sim}{\longleftarrow} & \left({\mathcal C}{\mathcal C}_*^{tilt}(\Gamma,v)\right)_{Z_v} & \overset{\sim}{\longrightarrow} &  \left({\mathcal C}^{alt}_*(\Gamma,v)\right)_{Z_v}
\end{array}
$$
of ind-Fr\'echet-complexes.\\
{\bf Step 3:}\\
It is easily verified that the canonical chain map
$$
\begin{array}{ccc}
p_v:\,CC_*({\mathbb C}\langle\Gamma\rangle\rtimes U:{\mathbb C}\rtimes U)_{\{v\}} & \longrightarrow & CC_*({\mathbb C}\rtimes\Gamma)_{[v]} \\
 & & \\
 u_v\langle g_0\rangle d\langle g_1\rangle\ldots d\langle g_n\rangle & \mapsto & (g_n^{-1}vg_0)d(g_0^{-1}g_1)\ldots d(g_{n-1}^{-1}g_n) \\
 & & \\
 u_vd\langle g_1\rangle\ldots d\langle g_n\rangle & \mapsto & d(g_n^{-1}vg_1)d(g_1^{-1}g_2)\ldots d(g_{n-1}^{-1}g_n). \\
\end{array}
\eqno(5.12)
$$
induces isomorphisms 
$$
CC_*({\mathbb C}\langle\Gamma\rangle\rtimes U:{\mathbb C}\rtimes U)_{\{v\}}/Z_v\,\overset{\simeq}{\longrightarrow}\,
CC_*({\mathbb C}\rtimes\Gamma)_{[v]}
\eqno(5.13)
$$
of complexes of vector spaces. The norms $\parallel-\parallel_{(\lambda,N,k)}$ on $CC_*({\mathbb C}\langle\Gamma\rangle\rtimes U:{\mathbb C}\rtimes U)_{\{v\}}$ are weighted $\ell^1$-norms on the linear space spanned by the set $\Delta_\bullet^{cyc}(\Gamma)$ of cyclic simplices. The group $Z_v$ acts on this set and the $v$-weights of its elements are constant along the $Z_v$-orbits. This implies that the quotient norm on the space $CC_*({\mathbb C}\langle\Gamma\rangle\rtimes U:{\mathbb C}\rtimes U)_{\{v\}}/Z_v$ of $Z_v$-coinvariants is a weighted $\ell^1$-norm on the set of orbits, the weight of an orbit being equal to the weight of any of its elements. Thus the quotient norm on $CC_*({\mathbb C}\langle\Gamma\rangle\rtimes U:{\mathbb C}\rtimes U)_{\{v\}}/Z_v$ equals the norm $\parallel-\parallel_{(N,k)}$ on 
$CC_*(\ell^1_\lambda(\Gamma))_{[v]}$ under (5.12) and we obtain
an isomorphism
$$
CC_*({\mathbb C}\langle\Gamma\rangle\rtimes U:{\mathbb C}\rtimes U)_{(v,\lambda,N)}/Z_v\,\overset{\simeq}{\longrightarrow}\,
(CC_*(\ell^1_\lambda(\Gamma))_N)_{[v]},\,\,\,\lambda,N>1,
$$
of Fr\'echet-complexes and finally an isomorphism 
$$
{\mathcal C}{\mathcal C}'_*({\mathbb C}\langle\Gamma\rangle\rtimes U:{\mathbb C}\rtimes U)_{\{v\}}/Z_v\,\overset{\simeq}{\longrightarrow}\,
\underset{\lambda\to 1}{\lim}\,CC^\omega_*(\ell^1_\lambda(\Gamma))_{[v]}
\eqno(5.14)
$$
of the corresponding ind-Fr\'echet complexes.\\
{\bf Step 4:}\\
We claim that the complexes of coinvariants $C^{R}_*(\Gamma,v)^{alt}/Z_v$ are finite dimensional for all $R>0$. 
For word-hyperbolic groups this is clear as $C_*(\Gamma,v)=C_*(Z_v,{\mathbb C})$ by definition. For $\Gamma$ a discrete isometry group of a $CAT(0)$-space
we show, in the notations of Section 2,  that the set 
$$
\Delta_0(\Gamma,v)\,=\,\{g\in\Gamma,\,g\cdot D\cap X^v\neq\emptyset\}
$$
is the union of finitely many $Z_v$-orbits. Clearly $\Delta_0(\Gamma,v)$  is stable under left translation by $Z_v$ because the fixed point set $X^v$ is so. Its image under the map 
$$
\Gamma\,\longrightarrow\,Z_v\backslash\Gamma\simeq [v]\subset\Gamma,\,\,\,g\,\mapsto\,g^{-1}vg
$$
consists of elements of word-length at most $2\,diam(D)$ and is thus finite.
The ind-complex $ \left({\mathcal C}^{alt}_*(\Gamma,v)\right)_{Z_v}$ is therefore independent of the chosen norms on it. \\
{\bf Step 5:}\\
The previous argument also shows that the proper and isometric $Z_v$-action on the non-empty and complete $CAT(0)$-space $X^v$ is cocompact. 
It is well known that under these assumption the constant $Z_v$-module $\mathbb C$ admits a finite resolution $P_*$ by finitely generated, projective $Z_v$-modules (take the nerve of a $Z_v$-invariant Leray-cover of $X^v$, obtained from a finite Leray-cover of $X^v/Z_v$). As any two projective resolutions of the constant $Z_v$-module are chain-homotopy equivalent there is a chain homotopy equivalence 
$C_*^{alt}(\Gamma,v)\,\to\,P_*$ of complexes of $Z_v$-modules which restricts to a homotopy equivalence
$$
{\mathcal C}^{alt}_*(\Gamma,v)\,=\, "\underset{R\to\infty}{\lim} "C^{R}_*(\Gamma,v)^{alt}\,\overset{\sim}{\longrightarrow}\,P_*
$$
of ind-complexes of $\Gamma$-modules. Passing to coinvariants one obtains a chain-homotopy equivalence
$$
 \left({\mathcal C}^{alt}_*(\Gamma,v)\right)_{Z_v}\,\overset{\sim}{\longrightarrow}\,P_*/Z_v.
$$
The right hand side is a constant, finite dimensional ind-complex. It is chain-homotopy equvalent to the constant ind-complex with vanishing differentials which equals 
the direct sum of the homology groups of $P_*/Z_v$ (graded by the parity of degrees). Altogether we obtain 
$$
\underset{\lambda\to (1,\infty)}{\lim}\,CC^\omega_*(\ell^1_\lambda(\Gamma))_{[v]}\,\overset{\sim}{\longrightarrow}\,\left({\mathcal C}^{alt}_*(\Gamma,v)\right)_{Z_v}\,\overset{\sim}{\longrightarrow}
$$
$$
\overset{\sim}{\longrightarrow}\,H_*(P_*/Z_v)\,\overset{\simeq}{\longrightarrow}\,H_*(Z_v,{\mathbb C})\,\overset{\simeq}{\longrightarrow}\,H_*(\Gamma,Ind_{Z_v}^\Gamma{\mathbb C})\,\overset{\simeq}{\longrightarrow}\,H_*(\Gamma,{\mathbb C}[v]),
$$
where the isomorphism in the middle follows from Shapiro's lemma. In the case of a word-hyperbolic group one may take $P_*=C_*^R(Z_v,{\mathbb C})^{alt}$ 
for $R>>0$ large. 
\end{proof}

\subsection{The hyperbolic part}
Suppose that $\Gamma$ acts properly, isometrically and cocompactly on a $CAT(0)$-space. The notations of Section 2 will be understood.
Let $v\in\Gamma$ be an element of infinite order. Its stable length or minimal displacement 
$$
\ell_s(v)\,=\,\underset{y\in X}{\Inf}\,d_X(y,vy)
\eqno(5.15)
$$
is strictly positive. The stable length is invariant under conjugation and 
$$
\underset{\underset{\vert v\vert=\infty}{v\in\Gamma}}{\Inf}\,\ell_s(v)>0
\eqno(5.16)
$$
because $\Gamma$ acts cocompactly on $X$.
The set of points of minimal displacement
$$
Min(v)\,=\,\{x'\in X,\,d(x',vx')=\ell_s(v)\}
\eqno(5.17)
$$
is non-empty by the cocompactness of the action, convex and closed. It is the disjoint union of infinite geodesics which are stable under the action of $v$.
After the choice of a base point $y\in Min(v)$ one obtains a canonical isometry
$$
\Phi_y:Min(v)\,\simeq\,{\mathbb R}\times Y
$$
of the set of minimal displacement with the metric product of the real line and a complete $CAT(0)$-space $Y$ \cite{BH}. It is characterized by the following facts:
\begin{itemize}
\item If $\pi:\,Min(v)\,\to\,L_{y}$ is the ``orthogonal'' projection of $Min(v)$ onto the unique $v$-stable geodesic passing through $y$, then $\Phi_y(\pi^{-1}(y_0))=\{0\}\times Y$.\\
\item The $v$-stable geodesics of $Min(v)$ correspond under $\Phi_y$ to the fibres of the canonical projection of ${\mathbb R}\times Y$ onto $Y$.\\
\end{itemize}
The action of $v$ on $Min(v)$ corresponds to the action $(t,y)\mapsto (t+\ell_s(v),y)$ on ${\mathbb R}\times Y$.
We let $c_{(v,y)}:\Gamma\to{\mathbb R}$ be the map given by the composition
$$
c_{(v,y)}:\,\Gamma\,\overset{-\cdot x}{\longrightarrow}\,X\,\overset{\pi'_{Min(v)}}{\longrightarrow}\,Min(v)\,\overset{\Phi_y}{\longrightarrow}\,{\mathbb R}\times Y\,\overset{\pi_1}{\longrightarrow}\,{\mathbb R}
\,\overset{\cdot\ell_s(v)^{-1}}{\longrightarrow}\,{\mathbb R}.
\eqno(5.18)
$$
\\
Suppose now that $(\Gamma,S)$ is a $\delta$-hyperbolic group. Let $v$ be an element of infinite order which is of minimal word length in its conjugacy class. Fix a geodesic segment $\sigma$ joining the vertices $e$ and $v$ in the Cayley-graph ${\mathcal G}(\Gamma,S)$ and let $L\subset\vert{\mathcal G}(\Gamma,S)\vert$ be the infinite segment given by the union of all $U$-translates of $\sigma$. We equip $L$ with the maximal metric for which each edge is isometric to the unit interval and let $\sigma': L\to{\mathbb R}$ be the isometry which sends $e$ to 0 and $v$ to its word-length. Put
$$
c_{(v,L)}:\,\Gamma\,\overset{\pi_L}{\longrightarrow}\,C_0(L,{\mathbb C})\,\overset{\sigma'}{\longrightarrow}\,C_0({\mathbb R},{\mathbb C})\,
\,\overset{\cdot\ell_S(v)^{-1}}{\longrightarrow}\,C_0({\mathbb R},{\mathbb C}),
\eqno(5.19)
$$
\\
where $\ell_s(v)=\underset{n\to\infty}{\lim}\,\frac{\ell(v^n)}{n}>0$. The maps $c_{(v,y)}$ and $c_{(v,L)}$ intertwine the action of the infinite cyclic subgroup $U$ of $\Gamma$ generated by $v$ and the translation action of ${\mathbb Z}$ on the real line 
 and descend therefore to algebra homomorphisms
$$
\langle c_{(v,y)}\rangle,\,\langle c_{(v,L)}\rangle:\,{\mathbb C}\langle\Gamma\rangle\rtimes U\,\longrightarrow\,{\mathbb C}\langle{\mathbb R}\rangle\rtimes{\mathbb Z}.
$$
\begin{lemma}
Let the subgroup ${\mathbb Z}\subset{\mathbb R}$ act on the real line by translations. Let $\tau_1$ be the trace on ${\mathbb C}\rtimes{\mathbb Z}$ corresponding to the characteristic function of the conjugacy class $\{1\}\subset{\mathbb Z}$. Then the linear functional\\ $\chi: CC_1({\mathbb C}\langle{\mathbb R}\rangle\rtimes{\mathbb Z}:{\mathbb C}\rtimes {\mathbb Z})_{\{1\}}\, \longrightarrow\, {\mathbb C},$ given on one-forms by
$$
\begin{array}{cccc}
 & u_1d\langle t_1\rangle & \mapsto & 1, \\
 & & & \\
 & u_1\langle t_0\rangle d\langle t_1\rangle & \mapsto & t_1-t_0\\
\end{array}
\eqno(5.20)
$$
and vanishing on forms of higher degrees 
satisfies $b\chi=0$ and\\ $B\chi\,=\,(\epsilon\rtimes id)^*(\tau_1),$ where 
$\epsilon\rtimes id:{\mathbb C}\langle{\mathbb R}\rangle\rtimes {\mathbb Z}\,\longrightarrow\,{\mathbb C}\rtimes {\mathbb Z}$ is the augmentation.
\end{lemma}

\begin{proof}
$$
b\chi(u_1\langle t_0\rangle d\langle t_1\rangle d\langle t_2\rangle)\,=\,\chi(u_1\langle t_0\rangle \langle t_1\rangle d\langle t_2\rangle)\,-\,
\chi(u_1\langle t_0\rangle d(\langle t_1\rangle\langle t_2\rangle))\,+\,\chi(\langle t_2\rangle u_1 \langle t_0\rangle d\langle t_1\rangle)
$$

$$
=\,\chi(u_1\langle t_1\rangle d\langle t_2\rangle)\,-\,
\chi(u_1\langle t_0\rangle d\langle t_2\rangle)\,+\,\chi(u_1\langle t_0\rangle d\langle t_1\rangle)
$$
$$
=\,(t_2-t_1)-(t_2-t_0)+(t_1-t_0)\,=\,0,
$$

$$
b\chi(u_1d\langle t_1\rangle d\langle t_2\rangle)\,=\,\chi(u_1\langle t_1\rangle d\langle t_2\rangle)\,-\,
\chi(u_1d(\langle t_1\rangle\langle t_2\rangle))\,+\,\chi(\langle t_2\rangle u_1 d\langle t_1\rangle)
$$

$$
=\,\chi(u_1\langle t_1\rangle d\langle t_2\rangle)\,-\,
\chi(u_1d\langle t_2\rangle)\,+\,\chi(u_1\langle t_2-1\rangle d\langle t_1\rangle)
$$

$$
=\,(t_2-t_1)-1+(t_1-(t_2-1))\,=\,0,
$$

$$
B\chi(u_1\langle t_0\rangle)\,=\,\chi(u_1d\langle t_0\rangle)\,=\,1\,=\,\tau_1((\epsilon\rtimes id)(u_1\langle t_0\rangle)).
$$
\end{proof}

\begin{lemma}
Let $v\in\Gamma$ be an element of infinite order, let $[v]$ be its conjugacy class, and let $U$ be the infinite cyclic group generated by $v$.
\begin{itemize}
\item[a)] 
Suppose that $\Gamma$ is a discrete isometry group of a $CAT(0)$-space as in \\Section 2 and fix a base point $y\in Min(v)\subset X$.
The cochain
$$
\chi_v\,=\,\chi_{(v,y)}\,=\,\langle c_{(v,y)}\rangle^*(\chi)\,\in\,CC^1({\mathbb C}\langle\Gamma\rangle\rtimes U:{\mathbb C}\rtimes U)_{\{v\}}
\eqno(5.21)
$$
is $Z_v$-invariant.
\item[b)]
Suppose that $\Gamma$ is word-hyperbolic and that $v$ is of minimal word length in its conjugacy class. Let $L$ be an infinite quasi-geodesic segment as constructed above.
Then the cochain
$$
\chi_v\,=\,\chi_{(v,L)}\,=\,\frac{1}{\vert N_v\vert}\,\underset{\overline{g}\in N_v}{\sum}\,g^*(\langle c_{(v,L)}\rangle^*(\chi))
\eqno(5.22)
$$
is $Z_v$-invariant. Here the average is taken over the\\ finite group $N_v=Z_v/U$.
\item[c)]
The $Z_v$-invariant cochains of a) and b) descend via (5.12) to cochains
$$
\overline{\chi}_{v}\,\in\,CC^1({\mathbb C}\rtimes\Gamma)_{[v]}.
\eqno(5.23)
$$
Their coboundary equals the trace on the group algebra corresponding to the characteristic function of the conjugacy class $[v]$.
\end{itemize}
\end{lemma}

\begin{proof}
The non-empty, closed and convex set $Min(v)$ of points of minimal displacement under $v$ is invariant under the action of $Z_v$. Therefore the ``orthogonal'' projection of $X$ onto $Min(v)$ commutes with the $Z_v$-action.
Under the identification $\Phi_y:Min(v)\overset{\simeq}{\longrightarrow}{\mathbb R}\times Y$ the $Z_v$-action on $Min(v)$ corresponds to a product action on ${\mathbb R}\times Y$ (see \cite{BH},II,6.8), where the action on the first factor is given by orientation preserving isometries, i.e. by translations. This implies that the cochain $\overline{\chi}_{(v,y)}$ is $Z_v$-invariant. Suppose now that $\Gamma$ is word-hyperbolic. It is well known that the infinite cyclic subgroup $U$ generated by $v$ is of finite index in its centralizer \cite{BH},III,3.10. The cochain $\langle c_{(v,L)\rangle}^*(\chi)$ is $U$-invariant, so that one may average it over the finite group 
$N_v=Z_v/U$ and arrives at a $Z_v$-invariant functional. In both cases the $Z_v$-invariant cochain $\chi_v$ descends, via (5.12),  to a cochain on the group algebra ${\mathbb C}\rtimes\Gamma$ as claimed.
By Lemma 5.5 its coboundary equals the trace on the group ring associated to the characteristic function of the conjugacy class $[v]$.
\end{proof}

\begin{lemma}
\begin{itemize}
The notations of Lemma 5.6 are understood.
\item[a)] Let $\Gamma$ be an isometry group of a $CAT(0)$-space as in Section 2. Then
$$
\vert\chi_{v}(\xi)\vert\,\leq\,\frac{C(\lambda,N)}{\ell_s(v)}\cdot\parallel\xi\parallel_{(\lambda,N,k)}
\eqno(5.24)
$$
for all $\lambda,N>1,k\in{\mathbb N}$ and $\xi\in CC_1({\mathbb C}\langle\Gamma\rangle\rtimes U:{\mathbb C}\rtimes U)_{\{v\}}$.
\item[b)] Let $\Gamma$ be word-hyperbolic. Then
$$
\vert\chi_{v}(\xi')\vert\,\leq\,C(\Gamma,S,\delta,\lambda,N)\cdot\parallel\xi'\parallel_{(\lambda,N,k)}
\eqno(5.25)
$$
for all $\lambda,N>1,k\in{\mathbb N}$ and $\xi'\in CC_1({\mathbb C}\langle\Gamma\rangle\rtimes U:{\mathbb C}\rtimes U)_{\{v\}}$.
\end{itemize}
\end{lemma}

\begin{proof}
a) For $\langle g_0\rangle d\langle g_1\rangle\in\Delta_1^{cyc}(\Gamma)'$ we find 
$$
\,\vert\chi_y(\langle g_0\rangle d\langle g_1\rangle)\vert\,
=\,\frac{1}{\ell_s(v)}\vert\pi_1\circ\Phi_y\circ\pi'_{Min(v)}(g_1x)-\pi_1\circ\Phi_y\circ\pi'_{Min(v)}(g_0x))\vert
$$
$$
\leq\,\frac{1}{\ell_s(v)}d(\pi'_{Min(v)}(g_1x),\pi'_{Min(v)}(g_0x))
\,\leq\,\frac{1}{\ell_s(v)}d(g_1x,g_0x)\,\leq\,\frac{\vert\langle g_0\rangle d\langle g_1\rangle\vert_v}{\ell_s(v)}
$$
because orthogonal projections do not increase distances,
while for $d\langle g_1\rangle\in\Delta_1^{cyc}(\Gamma)''$ 
$$
\vert\chi_{v,y}(d\langle g_1\rangle)\vert\,=1\,
\leq\,\frac{d(g_1,vg_1)}{\ell_s(v)}\,=\,\frac{\vert d\langle g_1\rangle\vert_v}{\ell_s(v)}.
$$
Thus
$$
\vert\chi_{v,y}(\omega)\vert\,\leq\,C(\lambda,N)\frac{2^k}{N}\lambda^{\vert\omega\vert_v}\,=\,C(\lambda,N)\cdot\parallel\omega\parallel_{(\lambda,N,k)}
$$
for all cyclic simplices $\omega\in \Delta_*^{cyc}(\Gamma)$.\\
b) Recall \cite{BH},III H,1.13, that there exist constants $C_1(\delta),\,C_2(\delta)$ with the following properties.
If an element $v\in\Gamma$ of infinite order and minimal length in its conjugacy class satisfies $\ell(v)>C_1(\delta)$, then any geodesic segment joining 
$v^N$ and $v^{-N}$, $N>>0$, is at Hausdorff distance at most $C_2(\delta)$ from the corresponding segment in $L$. Consequently, in the notations of Section 2,
$$
d(g',h')\,\leq\,d(g,h)+C_3(\delta),\,\,\,\forall\,g'\in S'_L(g),\,h'\in S'_L(h),
$$
for some $C_3(\delta)>>0$ and all $g,h\in\Gamma$. A similar estimate holds for each indvidual conjugcy class by \cite{BH},III,3.9,3.11, and takes care of the finitely many conjugacy classes 
not covered by the previous argument. Therefore the reasoning of a) carries over (with modified constants) to hyperbolic groups.
\end{proof}

\begin{cor}
The trace $\tau_{hyp}$  on $\ell^1(\Gamma)$ associated to the characteristic function of the set of hyperbolic elements in $\Gamma$ is a coboundary in
$\underset{\underset{\lambda}{\leftarrow}}{\lim}\,CC_\omega^*(\ell^1_\lambda(\Gamma))$.
\end{cor}

\begin{proof}
The infimum of the stable lengths of non-torsion elements in $\Gamma$ (5.16) is strictly positive because the action of $\Gamma$ on $X$ is cocompact. 
As the bicomplex $CC_*(\ell^1_\lambda(\Gamma))_N$ is the topological direct sum of the bicomplexes $(CC_*(\ell^1_\lambda(\Gamma))_N)_{[v]}$, labeled by the conjugacy classes of $\Gamma$, it follows from 5.6 and 5.7 that
the cochain 
$$
\overline{\chi}\,=\,\underset{[v]}{\sum}\,\overline{\chi}_{v}
\eqno(5.26)
$$
obtained by summation of the cochains (5.21) or (5.22) over all conjugacy classes $[v]$ of hyperbolic elements in $\Gamma$, defines a bounded cochain on $CC^*(\ell^1_\lambda(\Gamma))_N$ for all $\lambda,N>1$ and thus on the ind-Fr\'echet complex $\underset{\lambda\to 1}{\lim}\,CC^*_\omega(\ell^1_\lambda(\Gamma))$ as well. Lemma 5.5 shows that its coboundary equals the trace $\tau_{hyp}$  associated to the characteristic function of the set of hyperbolic elements in $\Gamma$.
\end{proof}

\begin{prop}
The inhomogeneous part of the ind-Fr\'echet complex
$$
\underset{\lambda\to 1}{\lim}\, CC^\omega_*(\ell^1_\lambda(\Gamma))
\eqno(5.27)
$$
is contractible for the class of groups $\Gamma$ mentioned in Section 2.
\end{prop}

\begin{proof}
Recall the Eilenberg-Zilber Theorem in cyclic homology \cite{Pu3}, which states the existence of a natural chain-homotopy equivalence 
$$
\nabla: CC^\omega_*(A\otimes_\pi B)\,\overset{\sim}{\longrightarrow}\,CC^\omega_*(A)\otimes_\pi CC^\omega_*(B)
\eqno(5.28)
$$
for unital complex Banach-algebras $A,B$. If 
$\alpha\in CC^*_\omega(A)$ is a cochain in the topologically dual chain complex of $CC_*^\omega(A)$, then the slant product
$$
\backslash\alpha:\,CC^\omega_*(A\otimes_\pi B)\,\overset{\nabla}{\longrightarrow}\,CC^\omega_*(A)\otimes_\pi CC^\omega_*(B)
\,\overset{\alpha\otimes id}{\longrightarrow}\,CC^\omega_*(B)
\eqno(5.29)
$$
is a bounded linear map. It is a chain map iff $\alpha$ is a cocycle. If $\alpha=\tau$ happens to be a trace, then a chain map representaing the slant product is given by
$$
\begin{array}{cccc}
\backslash\tau: & CC^\omega_*(A\otimes_\pi B) & \longrightarrow & CC^\omega_*(B) \\
 & & & \\
 & (a_0\otimes b_0)d(a_1\otimes b_1)\ldots d(a_n\otimes b_n) & \mapsto & \tau(a_0a_1\ldots a_n)\cdot b_0db_1\ldots db_n.
\end{array}
\eqno(5.30)
$$
Consider the diagonal homomorphism
$$
\Delta:\,\ell^1_{\lambda^2}(\Gamma)\,\longrightarrow\,\ell^1_{\lambda}(\Gamma)\otimes_\pi \ell^1_{\lambda}(\Gamma)
$$
and the induced chain map 
$$
\Delta_*:\,\underset{\lambda\to 1}{\lim}\,CC^\omega_*(\ell^1_\lambda(\Gamma))\,\longrightarrow\,
\underset{\lambda\to 1}{\lim}\, CC^\omega_*(\ell^1_\lambda(\Gamma)\otimes_\pi\ell^1_\lambda(\Gamma))
$$ 
of ind Fr\'echet complexes. The ``cap-product''
$$
\cap_{\overline{\chi}}:\,
\underset{\lambda\to 1}{\lim}\, CC^\omega_*(\ell^1_\lambda(\Gamma))\,\overset{\Delta_*}{\longrightarrow}\,
\underset{\lambda\to 1}{\lim}\, CC^\omega_*(\ell^1_\lambda(\Gamma)\otimes_\pi\ell^1_\lambda(\Gamma))\,\overset{\backslash\overline{\chi}}{\longrightarrow}\,
\underset{\lambda\to 1}{\lim}\, CC^\omega_{*+1}(\ell^1_\lambda(\Gamma))
\eqno(5.31)
$$
is then a bounded linear map of ind-Fr\'echet-complexes. Its coboundary equals 
$$
\delta(\cap_{\overline{\chi}})\,=\,\cap_{\delta(\overline{\chi})}\,=\,\cap_{\tau_{hyp}}:\,
\underset{\lambda\to 1}{\lim}\, CC^\omega_*(\ell^1_\lambda(\Gamma))\,\to\,\underset{\lambda\to 1}{\lim}\, CC^\omega_*(\ell^1_\lambda(\Gamma)).
$$
As
$$
\cap_{\tau_{hyp}}:\underset{\lambda\to 1}{\lim} CC^\omega_*(\ell^1_\lambda(\Gamma))\,\longrightarrow\,\underset{\lambda\to 1}{\lim} CC^\omega_*(\ell^1_\lambda(\Gamma))
$$
equals the canonical projection onto $\underset{\lambda\to 1}{\lim}\, CC^\omega_*(\ell^1_\lambda(\Gamma))_{inhom}$ by (5.30), the latter is contractible.
\end{proof}

\subsection{Conclusion}

\begin{theorem} 
Let $\Gamma$ be a discrete group acting properly, isometrically and cocompactly on a $CAT(0)$-space or suppose that $\Gamma$ is word-hyperbolic. 
Let $\Gamma_{tors}\subset\Gamma$ be the subset of elements of finite order, equipped with the adjoint $\Gamma$-action. There is an isomorphism
$$
\underset{\lambda\to 1}{\lim}\,CC_*^\omega(\ell^1_\lambda(\Gamma))\,\overset{\simeq}{\longrightarrow}\,H_*(\Gamma,{\mathbb C}\Gamma_{tors}) 
\eqno(5.32)
$$
in the chain-homotopy category of ind-Fr\'echet-complexes, where the right hand side is viewed as constant, finite dimensional ind-complex with vanishing differentials.
\end{theorem}

\begin{proof}
This follows  from Prop 5.4, Prop 5.9, and the fact that the groups studied here contain only finitely many conjugacy classes of torsion elements \cite{BH},III,1.1,3.2. 
\end{proof}

\begin{theorem} 
Let $\Gamma$ be a discrete group acting properly, isometrically and cocompactly on a $CAT(0)$-space or suppose that $\Gamma$ is word-hyperbolic. 
Let ${\mathcal C}{\mathcal C}_*(\ell^1(\Gamma))$ be the analytic cyclic bicomplex of the group Banach algebra $\ell^1(\Gamma)$. Then there is an isomorphism
$$
{\mathcal C}{\mathcal C}_*(\ell^1(\Gamma))\,\overset{\simeq}{\longrightarrow}\,H_*(\Gamma,{\mathbb C}\Gamma_{tors}) 
\eqno(5.33)
$$
in the derived ind-category.
\end{theorem}
\begin{proof}
The canonical homomorphisms $\ell^1_\lambda(\Gamma)\to\ell^1_{\lambda'}(\Gamma)\to\ell^1(\Gamma)$ are compact for $1<\lambda'<\lambda$. Thus one obtains 
a canonical morphism
$$
\underset{\lambda\to 1}{\lim}\,CC_*^\omega(\ell^1_\lambda(\Gamma))\,\longrightarrow\,{\mathcal C}{\mathcal C}_*(\ell^1(\Gamma))
$$
in the chain-homotopy category of ind-Fr\'echet complexes. 
The Approximation Theorem \cite{Pu2}, 6.13, and in particular its proof show that it becomes an isomorphism in the derived ind-category. The assertion follows then from the previous theorem.
\end{proof}

\begin{theorem}
Let $\Gamma,\Gamma'$ be discrete groups satisfying the hypothesis of the previous theorems and let $A,B$ be separable Banach algebras. Then
$$
HC_*^{loc}(\ell^1(\Gamma)\otimes_\pi A,\ell^1(\Gamma')\otimes_\pi B)\,\simeq\,Hom(H_*(\Gamma,{\mathbb C}\Gamma_{tors}),H_*(\Gamma',{\mathbb C}\Gamma'_{tors}))\otimes HC_*^{loc}(A,B).
\eqno(5.34)
$$
\end{theorem}

\begin{proof}
By the previous theorem there is an isomorphism 
$$
\Psi:{\mathcal C}{\mathcal C}_*(\ell^1(\Gamma))\,\overset{\simeq}{\longrightarrow}\,H_*(\Gamma,{\mathbb C}\Gamma_{tors})
$$
in the derived ind-category. 
It follows that 
$$
\Psi\otimes id: {\mathcal C}{\mathcal C}_*(\ell^1(\Gamma))\otimes_\pi{\mathcal C}{\mathcal C}_*(A)\,\overset{\simeq}{\longrightarrow}\,H_*(\Gamma,{\mathbb C}\Gamma_{tors})\otimes{\mathcal C}{\mathcal C}_*(A)
\eqno(5.35)
$$ 
is an isomorphism in $ind({\mathcal D})$ as well (the tensor product of a weakly contractible ind-complex with an auxiliary one is weakly contractible again).  Composing with the Eilenberg-Zilber morphism 
$$
\Delta:\,{\mathcal C}{\mathcal C}_*(\ell^1(\Gamma)\otimes_\pi A)\,\overset{\sim}{\longrightarrow}\,{\mathcal C}{\mathcal C}_*(\ell^1(\Gamma))\otimes_\pi{\mathcal C}{\mathcal C}_*(A),
$$
which is a chain-homotopy equivalence, we arrive at an isomorphism 
$$
{\mathcal C}{\mathcal C}_*(\ell^1(\Gamma)\otimes_\pi A)\,\overset{\simeq}{\longrightarrow}\,H_*(\Gamma,{\mathbb C}\Gamma_{tors})\otimes{\mathcal C}{\mathcal C}_*(A)
\eqno(5.36)
$$
in the derived ind-category. The same argument applies to the couple $(\Gamma',B)$. The theorem follows then from the fact that the set of morphisms between two objects of $ind({\mathcal D})$ is a complex vector space and that the composition of morphisms is bilinear.
\end{proof}

\begin{theorem}
Let $\Gamma$ be a discrete group acting properly, isometrically and cocompactly on a $CAT(0)$ space or suppose that $\Gamma$ is word-hyperbolic. 
Then the Chern-character \cite{Lo} from topological $K$-theory to local cycllic homology induces an isomorphism
$$
\begin{array}{cccc}
ch: & K_*(\ell^1(\Gamma))\otimes_{\mathbb Z}{\mathbb C} & \overset{\simeq}{\longrightarrow} & HC_*^{loc}(\ell^1(\Gamma)).
\end{array}
\eqno(5.37)
$$
\end{theorem}

\begin{proof}
By the higher index theorem of Connes-Moscovici there is a commutative diagram of ``assembly maps''
$$
\begin{array}{ccccc}
 & K^{top}_*(\underline{E\Gamma})\otimes_{\mathbb Z} & \overset{\mu}{\longrightarrow} & K_*(\ell^1(\Gamma))\otimes_{\mathbb Z}{\mathbb C} & \\
  & & & &  \\
  ch_{top} & \downarrow && \downarrow & ch \\
  & & & & \\
  & H_*(\Gamma,{\mathbb C}\Gamma_{tors}) & \longrightarrow & HC_*^{loc}(\ell^1(\Gamma)) & \\
\end{array}
\eqno(5.38)
$$
The vertical arrow on the left is an isomorphism by a classical theorm of Atiyah and Hirzebruch (see \cite{Co}), the upper horizontal arrow is an isomorphism by Vincent Lafforgue's Thesis \cite{La}, 
and the lower horizontal arrow is an isomorphism by Theorem 5.11. It follows that the vertical arrow on the right is an isomorphism as well.
\end{proof}

\begin{remark}
The only (partial) calculation of the local cyclic cohomology of a group Banach algebra prior to the present work was done in \cite{Pu4}, where its homogeneous part was determined for word-hyperbolic groups 
and used to verify the Kadison-Kaplansky idempotent conjecture for these groups. I want to point out that this paper contains an error. Proposition 2.9 is wrong so that a modified proof for Corollary 2.10 is needed.
It goes as follows. 
Replace the equivariant morphism of resolutions $\Phi':C_*(\Gamma,{\mathbb C})\to C_*^R(\Gamma,{\mathbb C})$ in Lemma 2.5.a) by the equivariant 
morphism $\varphi:C_*(\Gamma,{\mathbb C})\to C_*(X,{\mathbb C})$ (for $X=\vert\Delta_\bullet^R(\Gamma)\vert$) of \cite{Mi}, Proposition 12. Then corollary 2.10 holds
 for given $\lambda_1>1$ and $\lambda_0>1$ close enough to 1. The rest of the paper remains unchanged and all results 
 of chapters 3, 4 and 5 remain valid.
 \end{remark}

{\vskip8mm}
.\\
Institut de Mathématiques de Marseille (I2M),\\ 
UMR 7373 du CNRS, Campus de Luminy,\\ 
Universit\'e d'Aix-Marseille,\\ 
France


\begin{thebibliography}{Wittgenstein}  
\newcommand{\hotz}[1]{\bibitem[{#1}]{#1}} 

\hotz{BH} {\sc M.~Bridson, A.~Haefliger,} Metric spaces of non-positive\\ curvature, 
Springer Grundlehren 316, (1999)

\hotz{Co} {\sc A.~Connes,} Noncommutative Geometry, Academic Press (1994)

\hotz{CM} {\sc A.~Connes, H.~Moscovici,} Cyclic cohomology, the Novikov conjecture, and
 hyperbolic groups, Topology 29 (1990), 345-388

\hotz{CMR} {\sc J.~Cuntz, R.~Meyer, J.~Rosenberg} Topological and bivariant K-theory, Oberwolfach Seminars 36, Birkhaeuser (2007)

\hotz{KS} {\sc M.~Kashiwara, P.~Shapira,} Sheaves on manifolds,\\ Springer Grundlehren 292 (1990)

\hotz{La}  {\sc V.~Lafforgue,}  $K$-th\'eorie bivariante pour les alg\`ebres de Banach et conjecture de Baum-Connes, Invent. Math. 149 (2002), 1-95

\hotz{Lo}  {\sc J.L.~Loday,} Cyclic Homology, Springer Grundlehren 301 (1992)

\hotz{Mi}  {\sc I.~Mineyev,} Straightening and bounded cohomology of hyperbolic groups, Geom. Funct. Anal. 11 (2001), 807-839

\hotz{Ni} {\sc V.~Nistor,} Group cohomology and the cyclic cohomology of crossed products, Invent. Math. 99 (1990), 411-424

\hotz{Pu1} {\sc M.~Puschnigg,} Periodic cyclic homology of crossed products,\\ Proc. Symp. Pure Math. 105 (2023), 435-455

\hotz{Pu2} {\sc M.~Puschnigg,} Diffeotopy functors of ind-algebras and local cyclic cohomology, Docum. Math. 8 (2003), 143-245

\hotz{Pu3} {\sc M.~Puschnigg,} Explicit Product structures in cyclic homology theories, K-Theory 15 (1998), 323-345

\hotz{Pu4} {\sc M.~Puschnigg,} The Kadison-Kaplansky conjecture for\\ word-hyperbolic groups, Invent. Math. 149 (2002), 153-194

\end{thebibliography}
\end{document}